\input amstex
\input epsf

\documentstyle{amsppt}

\NoBlackBoxes
\magnification1200

\pageheight{9 true in}
\pagewidth{6.5 true in}

\def\Li{\text{\rm Li}}
\def\Re{\text{\rm Re}}
\def\lcm{\text{\rm lcm}}

\topmatter
\title
Prime Number Races
\endtitle
\author
Andrew Granville and Greg Martin
\endauthor
\endtopmatter

There's nothing quite like a day at the races.... The quickening of
the pulse as the starter's pistol sounds, the thrill when your
favorite contestant speeds out into the lead (or the distress if
another contestant dashes out ahead of yours), and the
accompanying fear (or hope) that the leader might change. And what
if the race is a marathon? Maybe one of the contestants will be
far stronger than the others, taking the lead and running at the
head of the pack for the whole race. Or perhaps the race will be
more dramatic, with the lead changing again and again for as long
as one cares to watch.

Our race involves the odd prime numbers, divided into two teams
depending on the remainder when they are divided by 4:

Mod 4 Race, Team 3:\ \ \ 03, 07, 11, 19, 23, 31, 43, 47, 59, 67,
71, 79, 83, ...

Mod 4 Race, Team 1:\ \ \ 05, 13, 17, 29, 37, 41, 53, 61, 73, 89,
93, 97, ...

\noindent In this Mod 4 Race\footnote{The integers in the arithmetic
progression $qn+a$ are often called the integers ``congruent to $a
\pmod q$'', and thus we have the ``Mod $q$ Race''.}, Team 3 contains
the primes of the
form $4n+3$, and Team 1 contains the primes of the form $4n+1$.
The Mod 4 Race has just two contestants and
    is quite some marathon  since it goes on forever!  From the data
above it appears that Team 3 is always in the lead; that is, up to
any given point, there seem to be at least as many primes of the
form $4n+3$ as there are primes of the form $4n+1$. Further data
seems to confirm our initial observations:
$$
\vbox{\offinterlineskip\hrule\halign{& \vrule#&
\strut\quad\hfil#\quad\cr height2pt&\omit&&\omit&\cr &\hfil&&\hfil
Number of primes \hfil && \hfil Number of primes
\hfil&\cr
    &$x$\hfil&&\hfil $4n+3$ up to $x$ \hfil && \hfil $4n+1$ up to $x$
\hfil&\cr height2pt&\omit&&\omit&\cr
\noalign{\hrule}height2pt&\omit&&\omit&\cr & \hfil $100$ \hfil &&
\hfil$13$ \hfil && \hfil $11$ \hfil&\cr\noalign{\hrule}\cr & \hfil
$200$ \hfil && \hfil$24$ \hfil && \hfil $21$
\hfil&\cr\noalign{\hrule}\cr & \hfil $300$ \hfil && \hfil$32$ \hfil
&& \hfil $29$ \hfil&\cr\noalign{\hrule}\cr & \hfil $400$ \hfil &&
\hfil$40$ \hfil && \hfil $37$ \hfil&\cr\noalign{\hrule}\cr & \hfil
$500$ \hfil && \hfil$50$ \hfil && \hfil $44$
\hfil&\cr\noalign{\hrule}\cr & \hfil $600$ \hfil && \hfil$57$ \hfil
&& \hfil $51$ \hfil&\cr\noalign{\hrule}\cr & \hfil $700$ \hfil &&
\hfil$65$ \hfil && \hfil $59$ \hfil&\cr\noalign{\hrule}\cr & \hfil
$800$ \hfil && \hfil$71$ \hfil && \hfil $67$
\hfil&\cr\noalign{\hrule}\cr & \hfil $900$ \hfil && \hfil$79$ \hfil
&& \hfil $74$ \hfil&\cr\noalign{\hrule}\cr & \hfil $1000$ \hfil &&
\hfil$87$ \hfil && \hfil $80$ \hfil&\cr\noalign{\hrule}\cr & \hfil
$2000$ \hfil && \hfil$155$ \hfil && \hfil $147$
\hfil&\cr\noalign{\hrule}\cr & \hfil $3000$ \hfil && \hfil$218$
\hfil && \hfil $211$ \hfil&\cr\noalign{\hrule}\cr & \hfil $4000$
\hfil && \hfil$280$ \hfil && \hfil $269$
\hfil&\cr\noalign{\hrule}\cr & \hfil $5000$ \hfil && \hfil$339$
\hfil && \hfil $329$ \hfil&\cr\noalign{\hrule}\cr & \hfil $6000$
\hfil && \hfil$399$ \hfil && \hfil $383$
\hfil&\cr\noalign{\hrule}\cr & \hfil $7000$ \hfil && \hfil$457$
\hfil && \hfil $442$ \hfil&\cr\noalign{\hrule}\cr & \hfil $8000$
\hfil && \hfil$507$ \hfil && \hfil $499$
\hfil&\cr\noalign{\hrule}\cr & \hfil $9000$ \hfil && \hfil$562$
\hfil && \hfil $554$ \hfil&\cr\noalign{\hrule}\cr & \hfil 10,000
\hfil && \hfil$619$ \hfil && \hfil $609$
\hfil&\cr\noalign{\hrule}\cr & \hfil 20,000 \hfil && \hfil$1136$
\hfil && \hfil $1125$ \hfil&\cr\noalign{\hrule}\cr & \hfil 50,000
\hfil && \hfil$2583$ \hfil && \hfil $2549$
\hfil&\cr\noalign{\hrule}\cr & \hfil 100,000 \hfil && \hfil$4808$
\hfil && \hfil $4783$ \hfil&\cr\noalign{\hrule}\cr
height2pt&\omit&&\omit&\cr}\hrule}
$$
\botcaption{Table 1} The number of primes of the form $4n+1$ and
$4n+3$ up to $x$.
\endcaption
\bigskip

Even in this extended data, the race remains close but Team 3
always seems to maintain a narrow lead. This phenomenon was first
observed in a letter written by Tch\'ebychev to M.~Fuss on 23rd
March 1853:

{\sl  \centerline{`` There is a notable difference in the
splitting of the prime numbers}

\centerline{ between the two forms $4n+3$, $4n+1$:}

    \centerline{ The first
form contains a lot more than the second.''}}
\smallskip
\noindent This bias is perhaps unexpected in light of an important
result in analytic number theory known as ``the prime number
theorem for arithmetic progressions''. This theorem tells us that
for any given modulus $q$, the primes tend to be equally split
amongst the various forms $qn+a$ where\footnote{Note that this
restriction is necessary, since every integer of the form $qn+a$
is divisible by the greatest common divisor of $a$ and $q$, hence
cannot be prime (except possibly for a single value $n$) if
gcd$(a,q)>1$.} $\gcd(a,q)=1$.  More precisely, we know that for
two such eligible values $a$ and $b$,
$$
\frac{\#\{\text{primes}\ qn+a \le x\}}{\#\{\text{primes}\ qn+b \le x\}}
\rightarrow 1\quad \text{as}\ x \rightarrow \infty.
\tag{1}
$$
This limit does not help us to predict who will win the
Mod $q$ Race\footnote{If the ratio $\#\{\text{primes}\ 4n+3 \le
x\}/\#\{\text{primes}\ 4n+1 \le x\}$ converged to a number greater
than 1 as $x\to\infty$, then we would know that $\#\{\text{primes}\
4n+3 \le x\}>\#\{\text{primes}\ 4n+1 \le x\}$ for all sufficiently
large $x$, so that in the long run Team 3 would always be ahead of Team
1. If it converged to a number less than 1, in the long run Team 1
would always be ahead of Team
3.}. In fact this {\sl asymptotic} result, the prime number
theorem for arithmetic progressions, does not inform us about any
of the fine details of these prime number counts, and so neither
verifies nor contradicts our observation that Team 3 seems to be
always ahead of Team 1.

It's time to come clean: we cheated in how we presented the data
in the table above. If we take the trouble to look at the counts
of Teams 1 and 3 for every value of $x$, not just the ones listed
in the table, we find that occasionally there is some drama to
this race, in that from time-to-time Team 1 takes the lead, but
then only briefly:

\noindent Team 1 takes the lead for the first time at
prime 26,861. However, since 26,863 is also prime, Team 3 draws equal
and then takes back the lead until

\noindent Team 1 gets ahead again at 616,841 and at various
numbers until 633,798.

\noindent Team 3 then takes back the lead until

\noindent Team 1 gets ahead again at 12,306,137 and at various
numbers until 12,382,326.

\noindent Team 3 then takes back the lead until

\noindent Team 1 gets ahead again at  951,784,481 and at various
numbers until 952,223,506.

    \noindent Team 3 then takes back the
lead until

\noindent Team 1 gets ahead again at  6,309,280,697 and at various
numbers until 6,403,150,362. \noindent Team 3 then takes back the
lead until

\noindent Team 1 gets ahead again at  18,465,126,217 and at
various numbers until 19,033,524,538.  \noindent Team 3 then takes
back the lead and holds on to it until at least twenty billion.

So there are, from time to time, more primes of the form $4n+1$
than of the form $4n+3$, but this lead is held only very briefly
and then relinquished for a long stretch. Nonetheless, given this
data, one might guess that $4n+1$ continues taking the lead from
time to time as we continue to watch this marathon. Indeed this is
the case:

\proclaim{J.E.\ Littlewood} {\rm (1914)}  There are arbitrarily
large values of $x$ for which there are more primes of the form
$4n+1$ up to $x$ than primes of the form $4n+3$. In fact there
    are arbitrarily large values of $x$ for which
$$
\#\{\text{primes}\ 4n+1 \le x\} - \#\{\text{primes}\ 4n+3 \le x\}
\ge \frac 12 \frac{\sqrt{x}}{\ln x} \ln \ln \ln x \tag{2}
$$
\endproclaim

At first sight, this seems to be the end of the story. However,
after seeing how infrequently Team 1 manages to hold on to the
lead, it is hard to put to rest the suspicion that Team 3 is in
the lead ``most of the time'', that there are more primes of the
form $4n+3$ up to $x$ than there are of the form $4n+1$, despite
Littlewood's result. In 1962, Knapowski and Tur\'an made a
conjecture that is consistent with Littlewood's result but also
bears out Tch\'ebychev's observation:

\proclaim{Conjecture}As $X \rightarrow \infty$, the percentage of
integers $x \le X$
for which there are more primes of the form $4n+3$ up to $x$ than of
the form $4n+1$ goes to 100\%.
\endproclaim

\noindent This conjecture may be paraphrased as

\centerline{``Tch\' ebychev was correct {\sl almost all} of the
time.''}

\noindent Looking at the above data in this light:
$$
\matrix X \ \text{in the}& \text{Maximum percentage of }\ x\leq X\
\text{for which there }\\
\underline{\text {Range}} & \underline{\text{are more primes}\ 4n+1\leq
x\ \text{than}\ 4n+3} \\
0-\text{26,860} & 0\% \\
0-\text{500,000} & \approx 0.01\% \\
0-10^7 & \approx 2.6\% \\
10^7-10^8 & \approx .6\% \\
10^8-10^9 & \approx .1\% \\
10^9-10^{10} & \approx 1.6\% \\
10^{10}-10^{11} & \approx 2.8\% \\
\endmatrix
$$
Does this persuade you that the Knapowski--Tur\'an conjecture is
likely to be true?  Is the percentage in the right column going to
$0$ as the values of $x$ get larger? The percentages are evidently
very low, but it is not obvious from this limited data that they
are tending towards 0. There was a wonderful {\sl Monthly} article
by Richard Guy some years ago entitled ``{\sl The Law of Small
Numbers}'', in which Guy pointed out several fascinating phenomena
that are ``evident'' for small integers yet disappear when one
examines bigger integers. Could this be one of those phenomena?

\subhead Another prime race: primes of the form $3n+2$ and
$3n+1$\endsubhead\smallskip

There are races between sequences of primes  other than that
between primes of the forms $4n+3$ and $4n+1$. For example, the
Mod 3 Race is between primes of the form $3n+2$ and primes of the
form $3n+1$. The race begins as follows:

Mod 3 Race, Team 2:\ \ 02, 05, 11, 17, 23, 29, 41, 47, 53, 59, 71,
83, 89, ...

Mod 3 Race, Team 1:\ \  07, 13, 19, 31, 37, 43, 61, 67, 73, 79,
97, ...

\noindent In this race Team 2 dashes out into an early lead which
it holds on to; that is, there seem to always be at least as
many primes of the form $3n+2$ up to $x$ as there are primes of
the form $3n+1$.  In fact Team 2 keeps in the lead up to ten
million and beyond\footnote{And this time we didn't cheat---Team 1 does
not get the
lead at some intermediate value of $x$ that we have not recorded
in the table.}.
$$
\vbox{\offinterlineskip\hrule\halign{& \vrule#&
\strut\quad\hfil#\quad\cr height2pt&\omit&&\omit&\cr
    &$x$\hfil&&\hfil Team 2 \hfil && \hfil
Team 1 \hfil&\cr height2pt&\omit&&\omit&\cr
\noalign{\hrule}height1pt&\omit&&\omit&\cr
\noalign{\hrule}height2pt&\omit&&\omit&\cr & \hfil $100$ \hfil &&
\hfil$13$ \hfil && \hfil $11$ \hfil&\cr\noalign{\hrule} \cr &
\hfil $200$ \hfil && \hfil$24$ \hfil && \hfil $21$
\hfil&\cr\noalign{\hrule}\cr & \hfil $300$ \hfil && \hfil$33$
\hfil && \hfil $28$ \hfil&\cr\noalign{\hrule}\cr & \hfil $400$
\hfil && \hfil$40$ \hfil && \hfil $37$
\hfil&\cr\noalign{\hrule}\cr & \hfil $500$ \hfil && \hfil$49$
\hfil && \hfil $45$ \hfil&\cr\noalign{\hrule}\cr & \hfil $600$
\hfil && \hfil$58$ \hfil && \hfil $50$
\hfil&\cr\noalign{\hrule}\cr & \hfil $700$ \hfil && \hfil$65$
\hfil && \hfil $59$ \hfil&\cr\noalign{\hrule}\cr & \hfil $800$
\hfil && \hfil$71$ \hfil && \hfil $67$
\hfil&\cr\noalign{\hrule}\cr & \hfil $900$ \hfil && \hfil$79$
\hfil && \hfil $74$ \hfil&\cr\noalign{\hrule}\cr & \hfil $1000$
\hfil && \hfil$87$ \hfil && \hfil $80$
\hfil&\cr\noalign{\hrule}\cr & \hfil $2000$ \hfil && \hfil$154$
\hfil && \hfil $148$ \hfil&\cr\noalign{\hrule}\cr & \hfil $3000$
\hfil && \hfil$222$ \hfil && \hfil $207$
\hfil&\cr\noalign{\hrule}\cr & \hfil $4000$ \hfil && \hfil$278$
\hfil && \hfil $271$ \hfil&\cr\noalign{\hrule}\cr & \hfil $5000$
\hfil && \hfil$338$ \hfil && \hfil $330$
\hfil&\cr\noalign{\hrule}\cr & \hfil $6000$ \hfil && \hfil$398$
\hfil && \hfil $384$ \hfil&\cr\noalign{\hrule}\cr & \hfil $7000$
\hfil && \hfil$455$ \hfil && \hfil $444$
\hfil&\cr\noalign{\hrule}\cr & \hfil $8000$ \hfil && \hfil$511$
\hfil && \hfil $495$ \hfil&\cr\noalign{\hrule}\cr & \hfil $9000$
\hfil && \hfil$564$ \hfil && \hfil $552$
\hfil&\cr\noalign{\hrule}\cr & \hfil $10000$ \hfil && \hfil$617$
\hfil && \hfil $611$ \hfil&\cr\noalign{\hrule}\cr & \hfil $20000$
\hfil && \hfil$1137$ \hfil && \hfil $1124$
\hfil&\cr\noalign{\hrule}\cr
height2pt&\omit&&\omit&\cr}\hrule}\qquad
\vbox{\offinterlineskip\hrule\halign{& \vrule#&
\strut\quad\hfil#\quad\cr height2pt&\omit&&\omit&\cr
    &$x$\hfil&&\hfil Team 2 \hfil && \hfil
Team 1 \hfil&\cr height2pt&\omit&&\omit&\cr
\noalign{\hrule}height1pt&\omit&&\omit&\cr
\noalign{\hrule}height2pt&\omit&&\omit&\cr & \hfil $30000$ \hfil
&& \hfil$1634$ \hfil && \hfil $1610$ \hfil&\cr\noalign{\hrule}\cr
& \hfil $40000$ \hfil && \hfil$2113$ \hfil && \hfil $2089$
\hfil&\cr\noalign{\hrule}\cr & \hfil $50000$ \hfil && \hfil$2576$
\hfil && \hfil $2556$ \hfil&\cr\noalign{\hrule}\cr & \hfil $60000$
\hfil && \hfil$3042$ \hfil && \hfil $3014$
\hfil&\cr\noalign{\hrule}\cr & \hfil $70000$ \hfil && \hfil$3491$
\hfil && \hfil $3443$ \hfil&\cr\noalign{\hrule}\cr & \hfil $80000$
\hfil && \hfil$3938$ \hfil && \hfil $3898$
\hfil&\cr\noalign{\hrule}\cr & \hfil $90000$ \hfil && \hfil$4374$
\hfil && \hfil $4338$ \hfil&\cr\noalign{\hrule}\cr & \hfil
$100000$ \hfil && \hfil$4807$ \hfil && \hfil $4784$
\hfil&\cr\noalign{\hrule}\cr & \hfil $200000$ \hfil && \hfil$8995$
\hfil && \hfil $8988$ \hfil&\cr\noalign{\hrule}\cr & \hfil
$300000$ \hfil && \hfil$13026$ \hfil && \hfil $12970$
\hfil&\cr\noalign{\hrule}\cr & \hfil $400000$ \hfil &&
\hfil$16967$ \hfil && \hfil $16892$ \hfil&\cr\noalign{\hrule}\cr &
\hfil $500000$ \hfil && \hfil$20804$ \hfil && \hfil $20733$
\hfil&\cr\noalign{\hrule}\cr & \hfil $600000$ \hfil &&
\hfil$24573$ \hfil && \hfil $24524$ \hfil&\cr\noalign{\hrule}\cr &
\hfil $700000$ \hfil && \hfil$28306$ \hfil && \hfil $28236$
\hfil&\cr\noalign{\hrule}\cr & \hfil $800000$ \hfil &&
\hfil$32032$ \hfil && \hfil $31918$ \hfil&\cr\noalign{\hrule}\cr &
\hfil $900000$ \hfil && \hfil$35676$ \hfil && \hfil $35597$
\hfil&\cr\noalign{\hrule}\cr & \hfil $1000000$ \hfil &&
\hfil$39266$ \hfil && \hfil $39231$ \hfil&\cr\noalign{\hrule}\cr &
\hfil $2000000$ \hfil && \hfil$74520$ \hfil && \hfil $74412$
\hfil&\cr\noalign{\hrule}\cr & \hfil $5000000$ \hfil &&
\hfil$174322$ \hfil && \hfil $174190$ \hfil&\cr\noalign{\hrule}\cr
& \hfil $10000000$ \hfil && \hfil$332384$ \hfil && \hfil $332194$
\hfil&\cr\noalign{\hrule}\cr height2pt&\omit&&\omit&\cr}\hrule}
$$
\botcaption{Table 2} The column labelled ``Team $j$'' contains the
number of primes of the form $3n+j$ up to $x$.
\endcaption
\bigskip

Perhaps our experience with the Mod 4 Race makes you skeptical
that Team 2 always dominates in this Mod 3 Race. If so, you are right
to be skeptical
because Littlewood's 1914 paper also applies to this race:
There are arbitrarily large values of $x$ for which there are more
primes of the form $3n+1$ up to $x$ than primes of the form
$3n+2$. (An inequality analogous to (2) also holds for this race.)
Therefore Team 1 takes the lead from Team 2 infinitely often, but
where is the first such value for $x$? We know it is beyond ten
million, and in fact Team 1 first takes the lead at
$$x = \text{608,981,813,029}.$$
This was discovered on Christmas day  1976 by Bays and Hudson. It
seems that Team 2 dominates in the Mod 3 Race even more than Team 3
dominates in the Mod 4 Race!

    \subhead Another prime race:  The last digit of a prime
\endsubhead\smallskip

Having heard that ``Team 2 retains the lead...'' for over a
half-a-trillion consecutive values of $x$, you might have
grown bored with your day at the races. So let's move on to a race
in which there are more competitors, perhaps making it less likely
that one will so dominate. One popular four-way race is between
the primes that end in 1, those that end in 3, those that end in 7
and  those that end in 9.
$$
    \matrix
\text{Last digit:}\hskip-6pt&\underline{\text{1}} &
\underline{\text{3}} &   \underline{\text{7}}
& \underline{\text{9}} &  \hskip9pt\text{Last digit:}\hskip-6pt  &
\underline{\text{1}} &  \underline{\text{3}} &  \underline{\text{7}} &
\underline{\text{9}}  \\
&\ &  3 &  7 & \ & \qquad & 101 &  103 & 107 & 109  \\
&11 &  13 &  17  & 19 & \qquad &  \ & 113  & \ & \\\
&\ &  23 & \  & 29 & \qquad & \ &  \ & 127 & \ \\
&31 &  \ &  37 & \ & \qquad & 131 & \ &  137 & 139\\
&41 &  43 &  47 & \ & \qquad & \ &  \ &  \ & 149 \\
&\ &  53 &  \ & 59 & \qquad & 151  & \  & 157 & \ \\
&61 &  \ &  67 & \ & \qquad & \  & 163 & 167 & \ \\
&71 &  73 &  \ & 79 & \qquad & \  & 173  & \ & 179 \\
&\ &  83 & \  & 89 &  & 181  & \  & \ & \ \\
&\ &  \ &  97 & \ & \qquad & 191  & 193  & 197 & 199  \\
\text{Up to 100:}\hskip-6pt&\overline{\ 5\ }  & \overline{\ 7\ } &
    \overline{\ 6\ } & \overline{\ 5\ } \
& \hskip9pt\text{Up to 200:}\hskip-6pt  & \overline{\ 10\ } &
\overline{\ 12\ }
& \overline{\ 12\ } & \overline{\ 10\ } \\
\endmatrix
$$
\botcaption{Table 3} The columns labelled ``Last digit $j$''
contain the primes  of the form $10n+j$ up to $100$, and  between
100 and 200, respectively.
\endcaption
\bigskip

On this limited evidence it seems that the two teams in the middle
lanes are usually ahead. However let us examine some more data
before making any bets:
$$
\matrix \underline{\ x\ } & \text{Last digit:} &
\underline{\text{1}}
& \underline{\text{3}} & \underline{\text{7}} & \underline{\text{9}} \\
100 & & 5 & 7 & 6 & 5 \\
200 & & 10 & 12 & 12 & 10 \\
500 & & 22 & 24 & 24 & 23 \\
1000 & & 40 & 42 & 46 & 38 \\
2000 & & 73 & 78 & 77 & 73 \\
5000 & & 163 & 172 & 169 & 163 \\
\text{10,000} & & 306 & 310 & 308 & 303 \\
\text{20,000} & & 563 & 569 & 569 & 559 \\
\text{50,000} & & 1274 & 1290 & 1288 & 1279 \\
\text{100,000} & & 2387 & 2402 & 2411 & 2390 \\
\text{200,000} & & 4478 & 4517 & 4503 & 4484 \\
\text{500,000} & & 10386 & 10382 & 10403 & 10365 \\
\text{1,000,000} & & 19617 & 19665 & 19621 & 19593 \\
\endmatrix
$$
\botcaption{Table 4} The column labelled ``Last digit $j$''
contains the   number of primes of the form $10n+j$ up to $x$.
\endcaption
\bigskip

    Except for a brief surge by Team 1 ahead of
Team 3 around $x=500,000$, it does appear that the two teams in
the middle lanes continue to share the lead.

After watching these races for a while, we begin to get the sense
that each race is somewhat predictable---not the tiny details,
perhaps, but the large-scale picture. If we were to watch a
similar race (the Mod $q$ Race for some number $q$ other than 4 or
3 or 10), we might expect that some of the teams would be much
stronger overall than others. But how could we predict which ones,
without watching for so long that everybody around us would
already be placing their bets on those same teams?

Before we can hope to understand  which team has the most primes
most often, we need to appreciate how those prime counts behave in
the first place. And before we can hope to understand the number
of primes up to $x$ of the form $qn+a$ (for given $a$ and $q$), we
should start by knowing how many  primes there are up to $x$. This
was perhaps the single most important question in 19th-century
number theory, a question which continued to engage researchers
throughout the 20th century and is still jealously guarding most
of its secrets from us today.

\head So what do we know about the count of primes up to $x$
anyway?\endhead

Although questions in number theory were not always mathematically
{\sl en vogue}, by the middle of the 19th century the problem of
counting primes had attracted the attention of well-respected
mathematicians such as Legendre, Tch\'ebychev, and the prodigious
Gauss. Although the best rigorous results of this time were due to
Tch\'ebychev, the prediction of Gauss   eventually led to a better
understanding of prime number races.

A query about the frequency with which primes occur elicited the
following response:
\smallskip
\centerline { \sl ``I  pondered this problem as a boy, and
determined that, at around $x$, the}

\centerline{{\sl  primes occur with density} $1/\ln x$.''
--- C.F.~Gauss (24th Dec 1849), letter to Encke}
\smallskip
\noindent This remark of Gauss can be interpreted as predicting
that
$$
\#\{\text{primes} \le x\} \approx \sum^{[x]}_{n=2} \frac{1}{\ln n}
\approx \int^x_2
\frac{dt}{\ln t} = \text{Li}(x),
$$
Let us compare Gauss's prediction\footnote{In fact, Euler made a
similar though less well-known ``prediction'' some years before
Gauss.} with the most recent count of the number of primes up to
$x$. (We shall wrote ``Overcount'' to denote $\Li(x)-\pi(x)$
rounded down to the nearest integer, which is the difference
between Gauss's prediction $\Li(x)$ and the function $\pi(x)$, the
true number of primes up to $x$.)

\centerline{\vbox{\offinterlineskip \halign{\vrule #&\ \ #\ \hfill
&& \hfill\vrule #&\ \ \hfill # \hfill\ \ \cr \noalign{\hrule} \cr
height5pt&\omit && \omit && \omit & \cr &$ \ \ x $&& $\pi(x) =
\#\{\text{primes} \le x\}$ &&  Overcount: $\Li(x)-\pi(x)$ & \cr
height5pt&\omit && \omit && \omit & \cr \noalign{\hrule} \cr
height3pt&\omit && \omit && \omit & \cr & $10^8$     &&  5761455
&& 753 & \cr height3pt&\omit && \omit && \omit & \cr & $10^9$ &&
50847534 && 1700  & \cr height3pt&\omit && \omit && \omit & \cr &
$10^{10}$  && 455052511 && 3103 & \cr height3pt&\omit && \omit &&
\omit & \cr &  $10^{11}$  && 4118054813 && 11587  & \cr
height3pt&\omit && \omit && \omit & \cr &  $10^{12}$  &&
37607912018 && 38262 & \cr height3pt&\omit && \omit && \omit & \cr
&  $10^{13}$  && 346065536839 && 108970 & \cr height3pt&\omit &&
\omit && \omit & \cr &  $10^{14}$  && 3204941750802 && 314889 &
\cr height3pt&\omit && \omit && \omit & \cr &  $10^{15}$  &&
29844570422669 && 1052618 & \cr height3pt&\omit && \omit && \omit
& \cr &  $10^{16}$  && 279238341033925 && 3214631 & \cr
height3pt&\omit && \omit && \omit & \cr &  $10^{17}$  &&
2623557157654233 && 7956588 & \cr height3pt&\omit && \omit &&
\omit & \cr &  $10^{18}$  && 24739954287740860 && 21949554 & \cr
height3pt&\omit && \omit && \omit & \cr &  $10^{19}$  &&
234057667276344607 && 99877774 & \cr height3pt&\omit && \omit &&
\omit & \cr &  $10^{20}$  && 2220819602560918840 && 222744643 &
\cr height3pt&\omit && \omit && \omit & \cr &  $10^{21}$ &&
21127269486018731928 && 597394253  & \cr height3pt&\omit && \omit
&& \omit & \cr &  $10^{22}$ && 201467286689315906290 && 1932355207
& \cr height3pt&\omit && \omit && \omit & \cr \noalign{\hrule}}} }
\botcaption{Table 5} Primes up to various $x$, and the overcount
in Gauss's prediction.
\endcaption
\bigskip

We can make several observations from these overcounts: First
notice that the width of the last column is always approximately
half the width of the middle column. In other words the overcount
seems to be about the square root of the main term $\pi(x)$
itself. Also the error term seems to always be positive, and so we
might guess that
$$
0 < \text{Li}(x) - \pi(x) < \sqrt{\pi(x)}
$$
for all $x$. In fact  the overcount appears to be monotonically
increasing which suggests that we might be able to better
approximate $\pi(x)$ by subtracting some smooth secondary term
from the approximating function $\Li(x)$, a question we shall
return to later. But, for now, let's  get back to the races\dots.

The function $\Li(x)$ does  not count primes, but it does seem to stay
close to
$\pi(x)$,  always remaining just in the lead. So for the race
between $\Li(x)$ and $\pi(x)$, we can ask: will $\Li(x)$
retain the lead forever? Littlewood's amazing result applies here
too:

\proclaim{Littlewood} {\rm (1914)}\ There are arbitrarily large values
of $x$ for which $\pi(x)>\Li(x)$, that is, for which
$$
\#\{\text{\rm primes}\le x\} > \int^x_2 \frac{dt}{\ln t}.
$$
\endproclaim

So what is the smallest $x_1$ for which $\pi(x_1) > \Li(x_1)$?
Skewes obtained an upper bound for $x_1$ from Littlewood's proof,
though not a particularly accessible bound. He proved:
$$
    \text{Skewes (1933):} \qquad x_1 < 10^{10^{10^{10^{34}}}},
$$
and to do even that he needed to make a significant assumption.
Skewes assumed the ``Riemann Hypothesis'', a conjecture that we
shall discuss a little later. For a long time, this ``Skewes
number'' was known as the largest number to which any
``interesting'' mathematical meaning could be ascribed. Skewes
later gave an upper bound that did not depend on any unproved
assumption, though at the cost  of making the numerical estimate
marginally more monstrous, and several improvements have been made
since then:
    $$ \align
    \text{Skewes (1955):}\qquad & x_1 < 10^{10^{10^{10^{1000}}}}\\
    \text{Lehman (1966):}\qquad & x_1 < 2 \times 10^{1165}\\
\text{te Riele (1987):}\qquad & x_1 < 6.658 \times 10^{370}\\
\text{Bays and Hudson (1999):}\qquad & x_1 < 1.3982 \times 10^{316}
\endalign
$$
As we shall discuss later, Bays and Hudson give compelling reasons
to believe that $x_1$ is actually  some integer close to $1.3982
\times 10^{316}$. This is an extraordinary claim! We can only
compute prime counts up to around $x=10^{22}$ with our current
technology and algorithms, and no significant improvements to this
are foreseen. So how can they make such a prediction of the enormous
value of $x_1$, when this prediction is so far beyond the
point where we can ever hope to compute $\pi(x)$ directly?

We shall explain this later. One  thought though: If the first
number $x$ for which $\pi(x)$ pulls ahead of $\Li(x)$ is so large
then surely the numbers $x$ for which $\pi(x)>\Li(x)$ are even
scarcer than the corresponding ``underdog'' numbers in the races
we examined earlier!

\subhead Accurately estimating the count of primes \endsubhead\smallskip

Up to the middle of the nineteenth century, every approach to
estimating $\pi(x) = \#\{\text{primes} \le x\}$ was relatively
direct, based upon elementary number theory and combinatorial
principles, or the theory of quadratic forms.
    In 1859, however, the great geometer Riemann took up the challenge of
counting
    the primes in a very different way. He wrote only one paper that 
could
be called
    number theory, but that one short memoir has had an impact lasting
nearly a century
    and a half already, and its ideas helped to define the subject we now
call analytic
    number theory.

Riemann's memoir described a surprising approach to the problem,
an approach using the theory of complex analysis, which was at that
time still very much a developing subject\footnote{Indeed,
Riemann's memoir on this number-theoretic problem was a
significant factor in the development of the theory of analytic
functions, notably their global aspects.}. This new approach of
Riemann's seemed to stray far away from the realm in which the
original problem lived. However, it had two key features:

\smallskip
\noindent $\bullet$ it was a potentially practical way to settle
the question once and for all;

\noindent $\bullet$ it made predictions that were similar, though
not identical, to the Gauss prediction. Indeed, as we shall see,
it suggested a secondary term to compensate somewhat for the
overcount we saw in the data above.

\smallskip
Riemann's method is too complicated to describe in its entirety
here,  but we shall take from it what we need to better understand
the prime count $\pi(x)$. To begin, we take the key prediction
from Riemann's memoir and restate it in entirely elementary
language:
$$
\lcm[1,2,3, \cdots, x] \approx e^x\ \text{as}\ x \rightarrow \infty.
$$
Now, one can easily verify  that
$$
\bigg( \prod_{p \le x} p \bigg) \times \bigg( \prod_{p^2 \le x} p
\bigg) \times \bigg( \prod_{p^3 \le x} p \bigg) \times \cdots =
\lcm[1,2,3, \cdots, x],
$$
since the power of any prime $p$ that divides the integer on
either  side of the equation is precisely the largest power of $p$
not exceeding $x$. Combining this identity with the preceding
approximation and taking natural logarithms of both sides, we
obtain
$$
\bigg( \sum_{p \le x} \ln p \bigg) + \bigg( \sum_{p^2 \le x}  \ln
p \bigg) + \bigg( \sum_{p^3 \le x} \ln p \bigg) + \cdots \approx
x\ \text{as}\ x \rightarrow \infty.
$$

Notice that the primes in the first sum are precisely the primes
counted by $\pi(x)$, the primes in the second sum are precisely
the primes counted by $\pi(x^{1/2})$, and so on. A technique
called partial summation allows us to ``take care of'' the $\ln p$
summand (which is analogous to using integration by parts to take
care of a $\ln x$ factor in an integrand). When partial summation
is applied to the above approximation, the result is
$$
\pi(x) + {\textstyle\frac12} \pi(x^{1/2}) + {\textstyle\frac13}
\pi(x^{1/3}) + \cdots \approx \int^x_2 \frac{dt}{\ln t} = \Li(x).
$$
If we ``solve for $\pi(x)$'' in an appropriate way, we find the
equivalent form
$$
\pi(x) \approx \Li(x) - {\textstyle\frac12} \Li(x^{1/2}) + \dots.
$$
Hence Riemann's method yields the same prediction as Gauss's
prediction, and it yields something extra as well---namely, it
predicts a secondary term that will hopefully compensate for the
overcount that we witnessed in the Gauss prediction. Let's review
the data and see how Riemann's prediction fares. ``Riemann's
overcount'' refers to $\Li(x)-\frac12\Li(\sqrt x)-\pi(x)$,
while ``Gauss's overcount'' refers to $\Li(x)-\pi(x)$ as before:
\smallskip

\centerline{\vbox{\offinterlineskip \halign{\vrule #&\ \ #\ \hfill
&& \hfill\vrule # &\ \ \hfill # \hfill\ \ \cr \noalign{\hrule} \cr
height5pt&\omit && \omit && \omit && \omit & \cr & $x$ && $\#
\{\text{primes}\le x$\} && Gauss's overcount && Riemann's
overcount & \cr height5pt&\omit && \omit && \omit && \omit & \cr
\noalign{\hrule} \cr height3pt&\omit &&\omit && \omit && \omit &
\cr &  $10^8$     &&  5761455&&    753&&131 & \cr height3pt&\omit
&& \omit && \omit && \omit & \cr & $10^9$     && 50847534&& 1700&&
-15  & \cr height3pt&\omit && \omit &&\omit && \omit & \cr &
$10^{10}$  && 455052511&& 3103&& -1711 & \cr height3pt&\omit &&
\omit && \omit && \omit & \cr &  $10^{11}$  && 4118054813&&
11587&& -2097 & \cr height3pt&\omit && \omit && \omit && \omit &
\cr &  $10^{12}$  && 37607912018&& 38262&& -1050 & \cr
height3pt&\omit && \omit && \omit && \omit & \cr &  $10^{13}$  &&
346065536839&& 108970&& -4944 & \cr height3pt&\omit && \omit &&
\omit && \omit & \cr &  $10^{14}$  && 3204941750802&& 314889&&
-17569 & \cr height3pt&\omit && \omit && \omit && \omit & \cr &
$10^{15}$  && 29844570422669&& 1052618&& 76456 & \cr
height3pt&\omit && \omit && \omit && \omit & \cr &  $10^{16}$  &&
279238341033925&& 3214631&& 333527 & \cr height3pt&\omit && \omit
&& \omit && \omit & \cr &  $10^{17}$  && 2623557157654233&&
7956588&& -585236 & \cr height3pt&\omit && \omit && \omit && \omit
& \cr &  $10^{18}$  && 24739954287740860&& 21949554&& -3475062 &
\cr height3pt&\omit && \omit && \omit && \omit & \cr &  $10^{19}$
&& 234057667276344607&& 99877774&& 23937697 & \cr height3pt&\omit
&& \omit && \omit && \omit & \cr &  $10^{20}$  &&
2220819602560918840&& 222744643&& -4783163 & \cr height3pt&\omit
&& \omit && \omit && \omit & \cr &  $10^{21}$  &&
21127269486018731928&& 597394253&& -86210244 & \cr height3pt&\omit
&& \omit && \omit && \omit & \cr &  $10^{22}$  &&
201467286689315906290&& 1932355207&& -126677992 & \cr
height5pt&\omit && \omit && \omit && \omit & \cr
\noalign{\hrule}}} } \botcaption{Table 6} Primes up to various
$x$, and Gauss's and Riemann's predictions.
\endcaption
\medskip

Riemann's prediction does seem to be a little better than that of
Gauss, although not a lot better. However, the fact that the error
in Riemann's prediction takes both positive and negative values
suggests that this might be about the best we can do.

Going back to the key prediction from Riemann's memoir, we can
calculate some data to test its accuracy:
$$
\matrix & \text{Nearest integer to} \\
\underline{x} & \underline{\ln(\lcm[1,2,3, \cdots, x])} &
\underline{\text{difference}}\\
100 & 94 & -6\\
1000 & 997 & -3\\
10000 & 10013 & 13\\
100000 & 100052 & 52\\
1000000 & 999587 & -413\\
\endmatrix
$$
Riemann's prediction has been made more precise over the years, and it
can now be expressed very explicitly as
$$
|\ln (\lcm[1,2, \cdots, x])-x| \leq 2 \sqrt{x} \ln^2 x\quad \text{for\
all}\ x \ge 100.
$$
In fact, this inequality is equivalent\footnote{As is
$|\Li(x)-\pi(x)|\leq \sqrt{x} \ln x$ for all $x\geq 3$.} to the
infamous {\sl Riemann Hypothesis}, perhaps the most prominent open
problem in mathematics. The Riemann Hypothesis is an assertion,
proposed by Riemann in his memoir and still unproven, about the
zeros of a certain function from complex analysis that is
intimately connected with the distribution of primes. To try to
give some idea of what Riemann did to connect the two seemingly
far-removed areas of number theory and complex analysis, we need
to talk about writing functions as combinations of ``waves''.

\subhead Doing the wave\endsubhead\smallskip

You have probably wondered what ``radio waves'' and ``sound
waves'' have to do with waves. Sounds don't usually seem to be
very ``wavy'', but rather are fractured, broken up, stopping and
starting and changing abruptly. So what's the connection? The idea
is that all sounds can be converted into a sum of waves. For
example, let's imagine that our ``sound'' is simply the gradually
ascending line $y=x-\frac12$, considered on the interval $0\le
x\le1$, which is shown in the first graph of Figure~1.

\midinsert {\centerline{ \epsfxsize=2.5truein\epsfbox{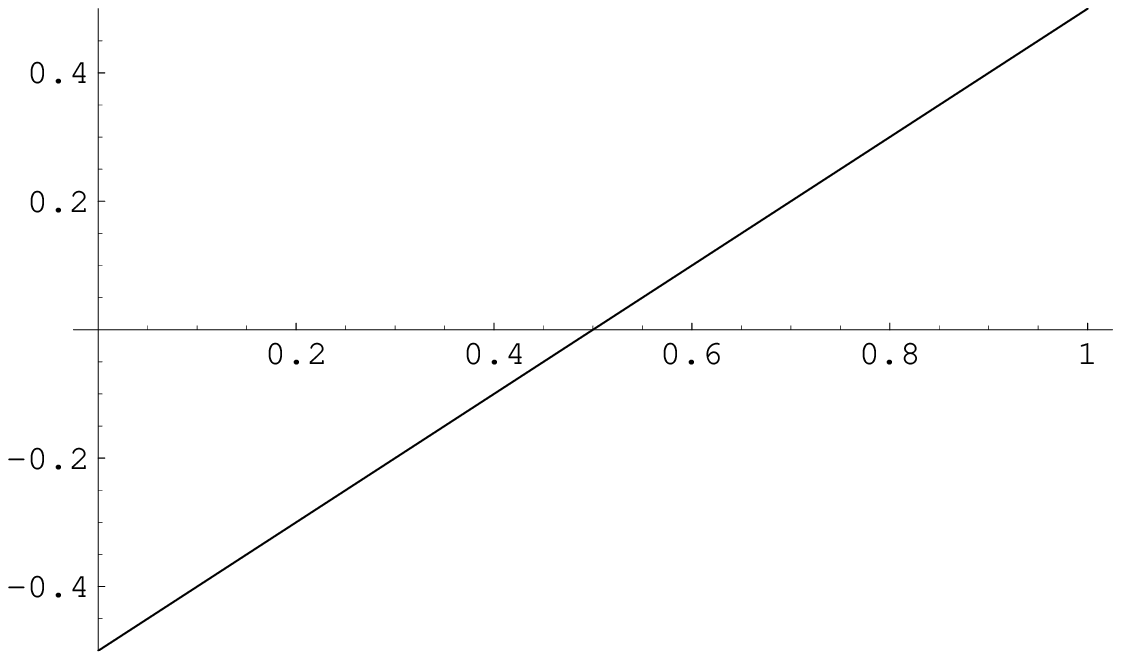}
\hskip.7in \epsfxsize=2.5truein\epsfbox{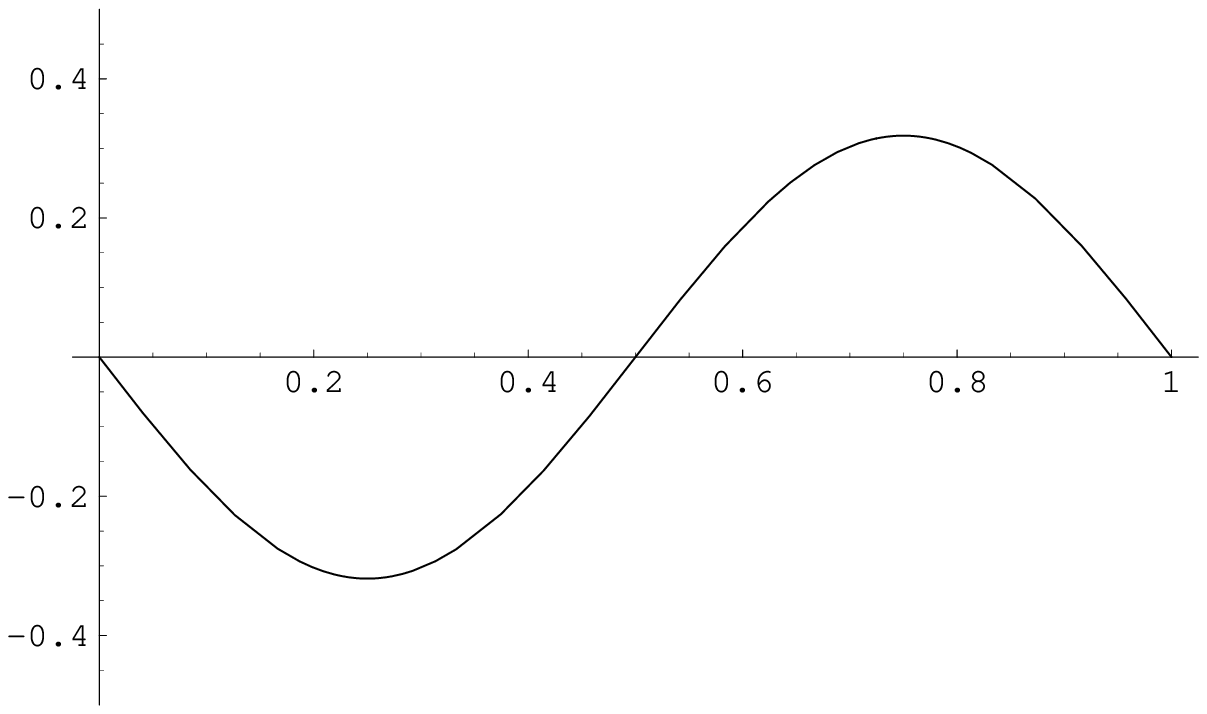} }}
\botcaption{Figure 1} The line $y=x-\frac12$ and the wave
$y=-\frac1{\pi}\sin 2\pi x$.
\endcaption
\endinsert

If we try to approximate this with a ``wave'', we can come pretty
close in the middle of the line using the function
$y=-\frac1{\pi}\sin 2\pi x$. However, as we see in the second
graph of Figure 1, the approximation is rather poor when
$x<\frac14$ or $x>\frac34$.

How can we improve this approximation? The idea is to ``add'' a
second wave to the first, this second wave going through two
complete cycles over the interval $0\le x\le1$ rather than only
one cycle. This corresponds to hearing the sound of the two waves
at the same time, superimposed; mathematically, we literally add
the two functions together. As it turns out, adding the function
$y=-\frac1{2\pi}\sin 4\pi x$ makes the approximation better for a
range of $x$-values which is quite a bit larger than the range of
good approximation we obtained with one wave, as we see in
Figure~2.

\midinsert {\centerline{ \epsfxsize=1.9truein\epsfbox{1wave.ps}
\hfill \hfill \raise1.16cm\hbox{+} \hfill
\epsfxsize=1.9truein\epsfbox{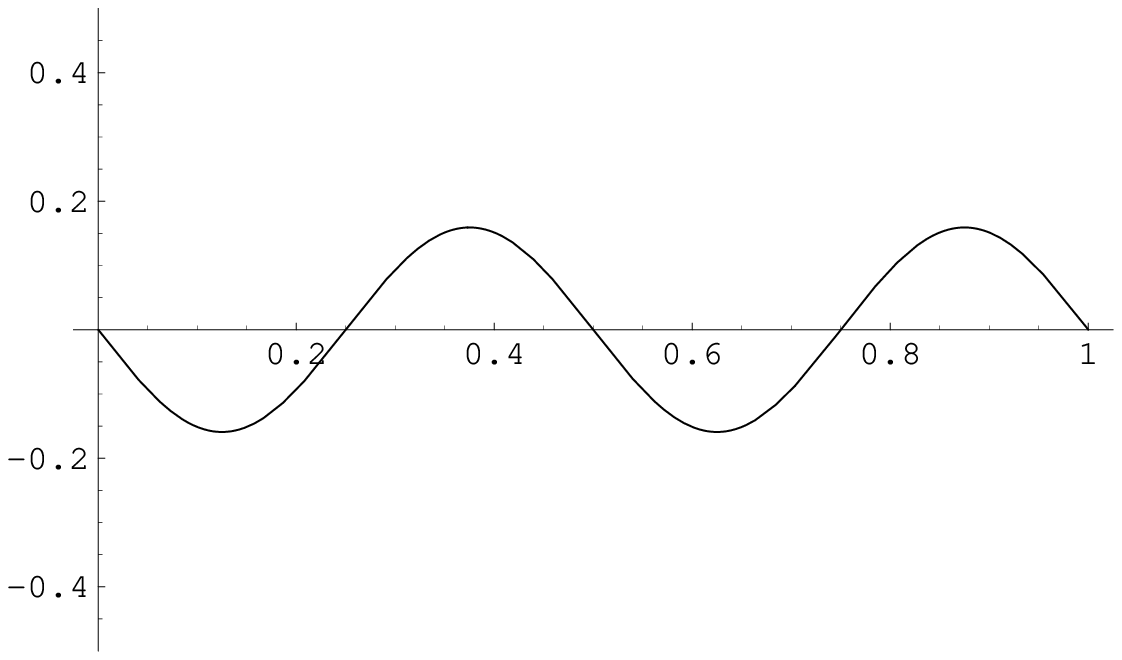} \hfill \hfill
\raise1.16cm\hbox{=} \hfill \epsfxsize=1.9truein\epsfbox{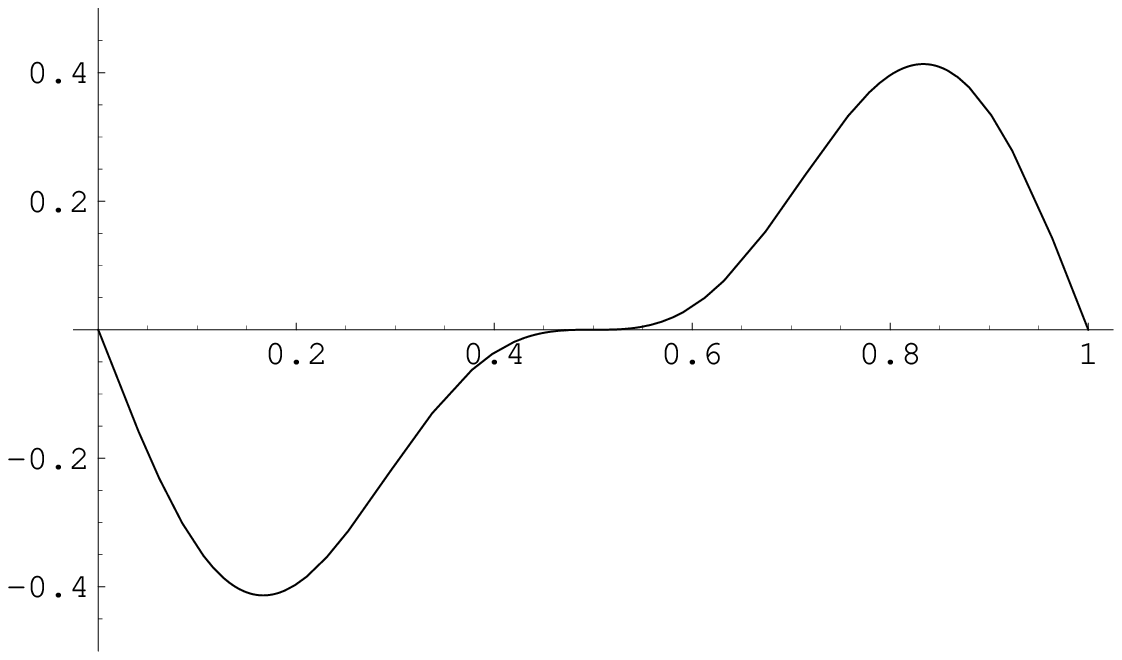} }}
\botcaption{Figure 2} Adding the wave $y=-\frac1{2\pi}\sin 4\pi x$ to
the wave
$y=-\frac1{\pi}\sin 2\pi x$
\endcaption
\endinsert

We can continue on like this, adding more and more waves that go
through 3, 4, or 5 complete cycles in the interval and so on, to get
increasingly better approximations to the original straight line.
The approximation we get by using one hundred superimposed waves
is really quite good, except near the endpoints 0 and 1.

\midinsert {\centerline{ \epsfxsize=6truein\epsfbox{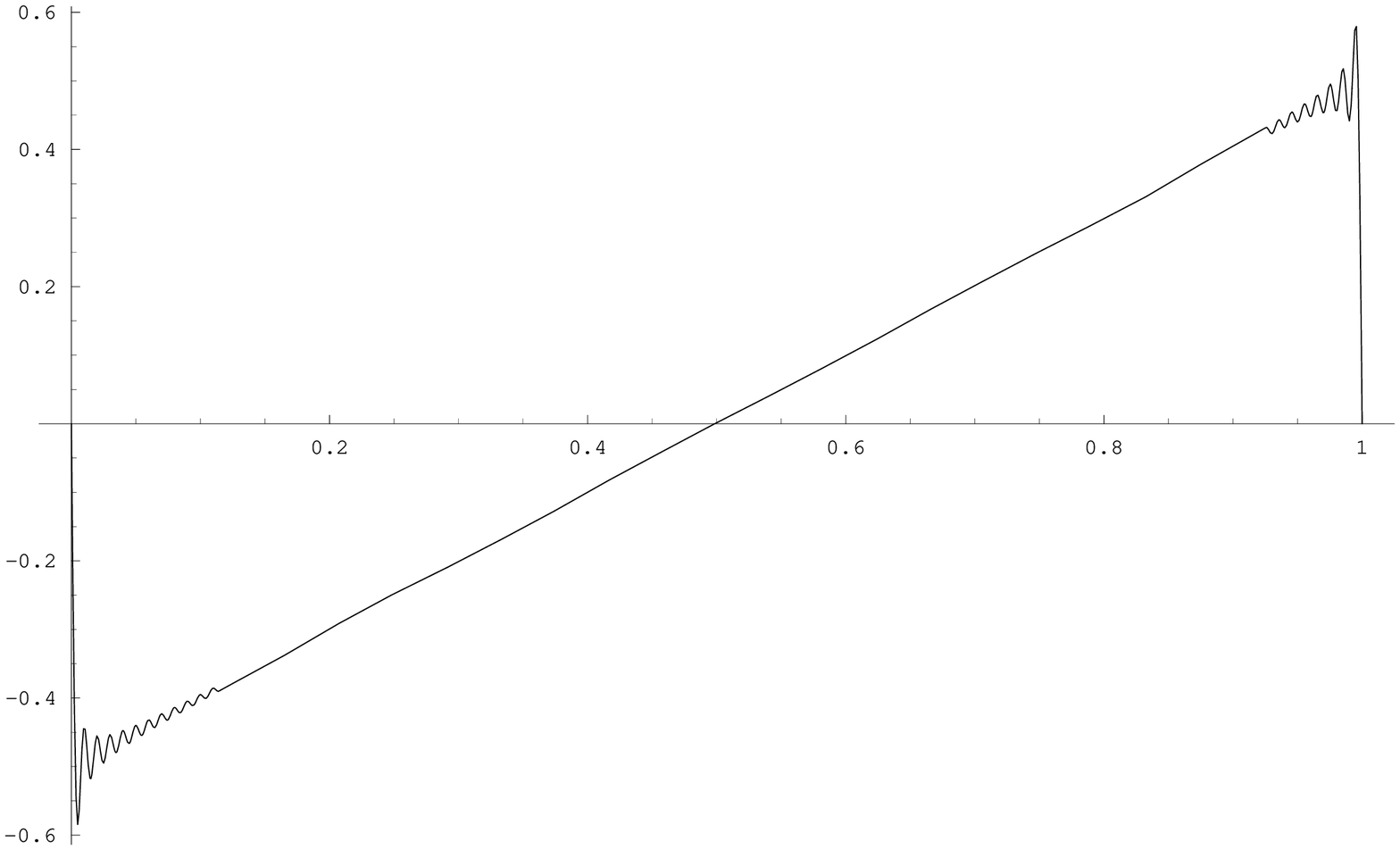}
}} \botcaption{Figure 3} The sum of one hundred carefully chosen
waves
\endcaption
\endinsert

If we were to watch these one hundred approximations being built
up one additional wave at a time, we would quickly be willing to
wager that the more waves we allowed ourselves, the better the
resulting approximation would be, perhaps becoming as accurate as
could ever be hoped for. As long as we allow a tiny bit of error
in the approximation (and shut our eyes to what happens very near
the endpoints), we can in fact construct a sufficiently close
approximation if we use enough waves. However, to get a
``perfect'' copy of the original straight-line, we would need to
use infinitely many sine waves---more precisely, the ones on the
right-hand side of the formula
$$
x-\frac12=-2\sum_{n=1,2,3,\dots} \frac{\sin(2\pi nx)}{2\pi n},
$$
which can be shown to hold for every $0<x<1$. (We can read  off,
in the $n=1$ and $n=2$ terms of this sum, the two waves that we
chose for Figures 1 and 2.) This formula is not of much practical
use, since we can't really transmit infinitely many waves at
once---but it's a gorgeous formula nonetheless!

In general, for any function $f(x)$ defined on the interval $0\le
x\le 1$ that is not ``too irregular'', we can find numbers $a_n$ and
$b_n$
such that $f(x)$ can be written as a sum of trigonometric
functions, namely
$$
f(x) = a_0 + \sum_{n=1,2,3,\dots} \big( a_n \cos (2 \pi nx) + b_n
\sin(2 \pi nx) \big).
$$
This formula, together with a way to calculate the coefficients
$a_n$ and  $b_n$, is one of the key identities from ``Fourier
analysis'', and it and its many generalizations are the subject of
the field of mathematics known as harmonic analysis. In terms of
waves, the numbers $2\pi n$ are the ``frequencies'' of the various
component waves (controlling how fast they go through their
cycles), while the coefficients $a_n$ and $b_n$ are their
``amplitudes'' (controlling how far up and down they go).

\subhead  Riemann's revolutionary formula\endsubhead\smallskip

Riemann's idea can be simply, albeit rather surprisingly, phrased
in the following terms:

\centerline {\sl Try counting the primes as a sum of waves.}

\noindent The precise formula he proposed is a bit too technical
for this article, but we can get a good sense of it from the
following approximation when $x$ is large. This formula, while
widely believed to be correct, has not yet been proved.
$$
\boxed{\frac{\int\limits_2^x{\frac{dt}{\ln t}} - \# \{\text{primes}\le
x\} }
{\sqrt x/\ln x} \approx 1 + 2\dsize\sum \Sb \text{all real numbers
} \gamma>0 \\ \text{such that } \frac 12 + i\gamma \\ \text{is a
zero of }\zeta(s)\endSb \ \ \frac{\sin(\gamma\ln x)}\gamma.}
\tag{3}
$$
The numerator of the left-hand side of this formula is the error
term when comparing the Gauss prediction $\Li(x)$ to the actual
count $\pi(x)$ for the number of primes up to $x$. We saw earlier
that the overcounts seemed to be roughly the size of the square
root of $x$, so the denominator $\sqrt x/\ln x$ appears to be an
appropriate thing to divide through by. The right side of the
formula bears much in common with our formula for $x-1/2$. It is a
sum of sine functions, with the numbers $\gamma$ employed in two
different ways in place of $2\pi n$: each $\gamma$ is used inside
the sine (as the ``frequency''), and the reciprocal of each
$\gamma$ forms the coefficient of the sine (as the ``amplitude'').
We even get the same factor of 2 in each formula. However, the
numbers $\gamma$ here are much more subtle than the
straightforward numbers $2\pi n$ in the corresponding formula for
$x-1/2$:

The {\sl Riemann zeta-function} $\zeta(s)$ is defined as
$$
\zeta(s)=\frac1{1^s}+\frac1{2^s}+\frac1{3^s}+\dots.
$$
Here $s$ is a complex number, which we shall write as $s=\sigma+it$
when we want to refer to its real and imaginary parts $\sigma$ and
$t$ separately. If $s$ were a real number, we know from first-year
calculus that the series in the definition of $\zeta(s)$ converges
if and only if $s>1$; that is, we can sum up the infinite series
and obtain a finite, unique value. In a similar way, it can be
shown that the series only converges for complex numbers $s$ such
that $\sigma>1$. But what about when $\sigma \leq 1$? How do we
get around the fact that the series does not sum up (that is,
converge)?

Fortunately, there is a beautiful phenomenon in the  theory of
functions of a complex variable called ``analytic continuation''.
It tells us that functions that are originally defined only for
certain complex numbers often have a unique ``sensible''
definition for other complex numbers. In this case, the definition
of $\zeta(s)$ given above only works when $\sigma>1$, but it turns
out that ``analytic continuation'' allows us to define $\zeta(s)$
for every  complex number $s$ other than $s=1$ (see [17] for
details).

This description of the process of analytic continuation looks
disconcertingly magical. Fortunately, there is a quite explicit
way to show how $\zeta(\sigma+it)$ can be ``sensibly'' defined for
the larger region $\sigma>0$, at least. We start with the
expression $(1-2^{1-s})\zeta(s)$ and perform some sleight-of-hand
manipulations:
$$
\align
    \left(1 - 2^{1-s}\right) \zeta(s) &= \Big( 1 - \frac{2}{2^s} \Big)
\zeta(s) = \zeta(s) - \frac{2}{2^s}\zeta(s) \\
&= \sum_{n\geq 1} \frac 1{n^s} - 2 \sum_{m\geq 1} \frac 1{(2m)^s} \\
&= \sum_{n\geq 1} \frac 1{n^s} - 2 \sum\Sb n\geq 1\\ n \
\text{even}\endSb \frac 1{n^s} =\  \sum\Sb m\geq 1\\ m \
\text{odd}\endSb \frac 1{m^s} -  \sum\Sb n\geq 1\\ n \ \text{even}\endSb
\frac 1{n^s} \\ &= \left(\frac{1}{1^s} - \frac{1}{2^s}\right) +
\left(\frac{1}{3^s} -
\frac{1}{4^s}\right) + \left(\frac{1}{5^s} - \frac{1}{6^s}\right) +
\cdots
\endalign
$$
Solving for $\zeta(s)$, we find that
$$
\zeta(s) = \frac 1 {(1 - 2^{1-s})} \left\{ \left(\frac{1}{1^s} -
\frac{1}{2^s}\right) + \left(\frac{1}{3^s} - \frac{1}{4^s}\right)
+ \left(\frac{1}{5^s} - \frac{1}{6^s}\right) + \cdots \right\}.
$$
All of these manipulations were valid for complex numbers
$s=\sigma+it$ for which $\sigma>1$. However, it turns out that the
infinite series in curly brackets actually converges whenever
$\sigma>0$. Therefore, we can take this last equation as the new
``sensible'' definition of the Riemann zeta-function on this
larger domain. Note that the special number $s=1$ causes the
otherwise innocuous factor of $1/(1 - 2^{1-s})$ to be undefined;
the Riemann zeta-function intrinsically has a problem there, one
that cannot be swept away with clever re-arrangements of the
infinite series.

    Riemann's formula (3) depends on the zeros of the analytic
continuation of $\zeta(s)$. The easiest zeros to identify are the
negative even integers; that is
$$
\zeta(-2)=0,\ \zeta(-4)=0,\ \zeta(-6)=0,\ \dots
$$
which are called the ``trivial zeros'' of the zeta-function. It
can be shown that any other complex zero $\sigma +it$ of the zeta
function (that is, a number satisfying $\zeta(\sigma+it)=0$) must
satisfy $0\le\sigma\le1$; these more mysterious zeros of the
zeta-function are called the ``non-trivial zeros''.

After some calculation Riemann observed that all of the
non-trivial zeros of the zeta-function seem to lie on the line
$\Re(s)=1/2$. In other words, he stated the remarkable hypothesis:
$$
\boxed{\text{If }\sigma+it \text{ is a complex number with }
0\le\sigma\le1 \text{ and } \zeta(\sigma+it)=0 \text{, then }
\sigma=\textstyle\frac12.}
$$
This assertion, which nobody has yet managed to prove, is the
famous {\sl Riemann Hypothesis}.

If the Riemann Hypothesis is indeed true, then we can write all
the non-trivial zeros of the zeta-function in the form $\rho =
\frac{1}{2} + i \gamma$ (together with their conjugates $\frac{1}{2}
- i \gamma$, since $\zeta(1/2+i\gamma)=0$ if and only if
$\zeta(1/2-i\gamma)=0$), where $\gamma$ is a positive real number.
These are
the mysterious numbers $\gamma$ appearing in the formula (3),
which holds if and only if the Riemann Hypothesis is true. There
is a similar formula if the Riemann Hypothesis is false, but it is
rather complicated and technically far less pleasant. The reason
is that the coefficients $1/\gamma$, which are constants in (3),
get replaced by functions of $x$. So we want the Riemann
Hypothesis to hold because it gives the simpler formula (3), and
that formula is a delight. Indeed (3) is similar enough to the
formulas for soundwaves for some experts to muse that (3) asserts
that ``{\sl the primes have music in them}''.

One might ask how we add up the infinite sum in (3)? Simple, add
up by order of ascending $\gamma$ values and it will work out.

Clever people have computed literally billions of zeros of
$\zeta(s)$, and every single zero that has been computed does
indeed sit squarely on the line $\sigma=1/2$. For example, the
non-trivial zeros closest to the real axis are $s=1/2+i\gamma_1$ and
$s=1/2-i\gamma_1$, where $\gamma_1\approx 14.1347 \dots$. We
believe that the positive numbers $\gamma$ occurring in the
non-trivial zeros look more or less like random real numbers, in
the sense that none of them are related to one another by simple
linear equations with integer coefficients (or even by more
complicated polynomial equations with algebraic numbers as
coefficients). However, since about all we know how to do is to
numerically approximate these non-trivial zeros to a given
accuracy, we cannot say much about the precise nature of the
numbers $\gamma$.

\subhead  Prime Race for $\pi(x)$ versus $\Li(x)$  \endsubhead

So, how do we use this variant on ``Fourier analysis'' to locate the
smallest $x$ for which
$$
\# \{ \text{primes} \le x \} > \int^x_2 \frac{dt}{\ln t} \ ?
$$
The idea is to approximate
$$
\frac{\Li(x) - \pi(x)}{\sqrt{x}/\ln x} = \frac{\int\limits^x_2
\frac{dt}{\ln t} - \#\{\text{primes} \le x\}}{\sqrt{x}/\ln x}
$$
by using the formula (3). Just as we saw when approximating the
line $x-1/2$, we expect to obtain a good approximation here simply
by summing the formula on the right side of (3) over the first few
zeros of $\zeta(s)$ (that is, the smallest hundred, or thousand,
or million values of $\gamma$ with $\zeta(1/2+i\gamma)=0$,
depending upon the level of accuracy that we want). In other
words,
$$
\frac{\int\limits_2^x{\frac{dt}{\ln t}} - \#\{\text{primes} \le
x\}}{\sqrt
x/\ln x} \approx 1 + 2 \sum \Sb \frac 12 + i\gamma\ \text{is a zero of}
\ \zeta(s)\\ 0<\gamma <T\endSb \ \ \frac{\sin(\gamma\ln x)}\gamma ,
$$
where we can choose $T$ to include however many zeros we want.
Figure 4 is a graph of the function\footnote{When $x$ is large,
$\frac12\Li(\sqrt x) \approx \sqrt x/\ln x$. We have preferred to use
the latter expression because it consists of more familiar functions.
However, when dealing with ``small'' values of $x$ as in these graphs,
using the function $\frac12\Li(\sqrt x)$ lets us more clearly display
the similarities we want to show.}
$(\Li(x)-\pi(x)) \big/
(\frac12\Li(\sqrt{x}))$ as $x$ runs from $10^4$ to $10^8$, together
with three graphs of approximations using the first 10, 100, and
1000 values of $\gamma$ respectively. We see that the
approximations get better the more zeros we take.

\midinsert {
\centerline{ \epsfxsize=5truein\epsfbox{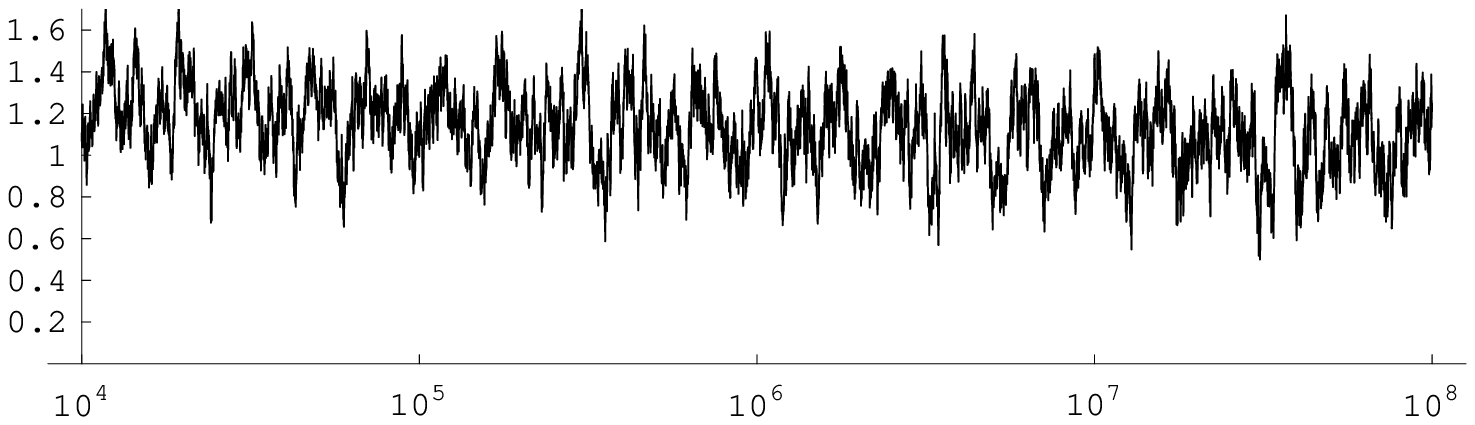} }
\vskip12pt
\centerline{ \epsfxsize=5truein\epsfbox{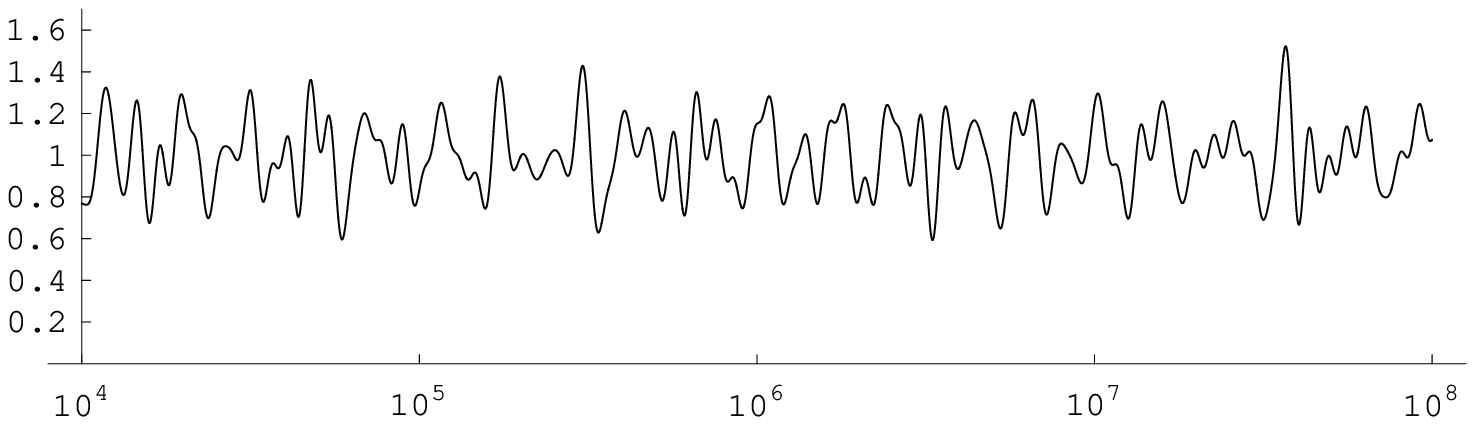} }
\centerline{ \epsfxsize=5truein\epsfbox{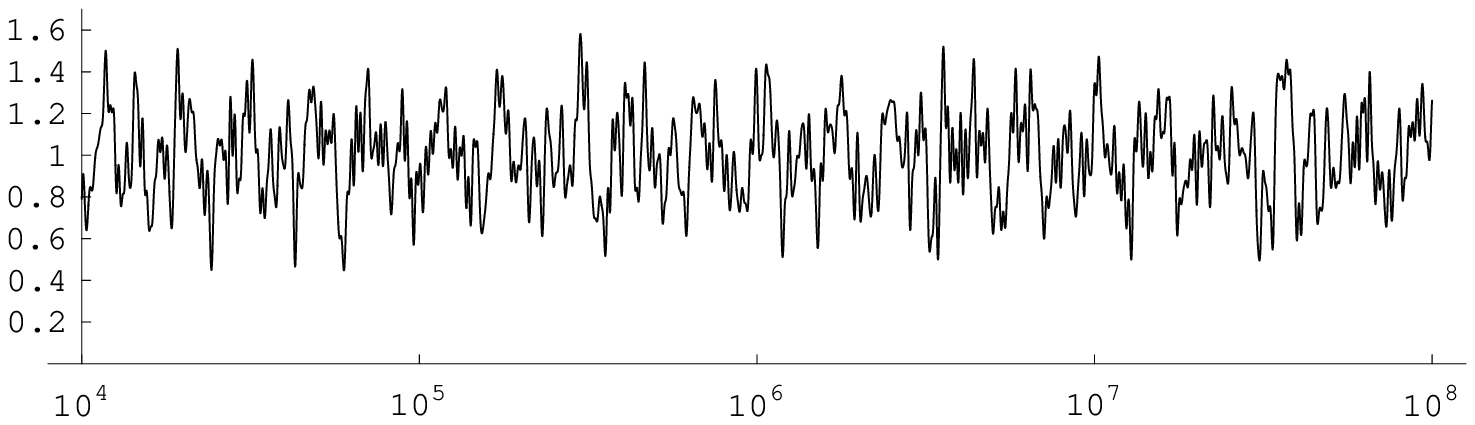} }
\centerline{ \epsfxsize=5truein\epsfbox{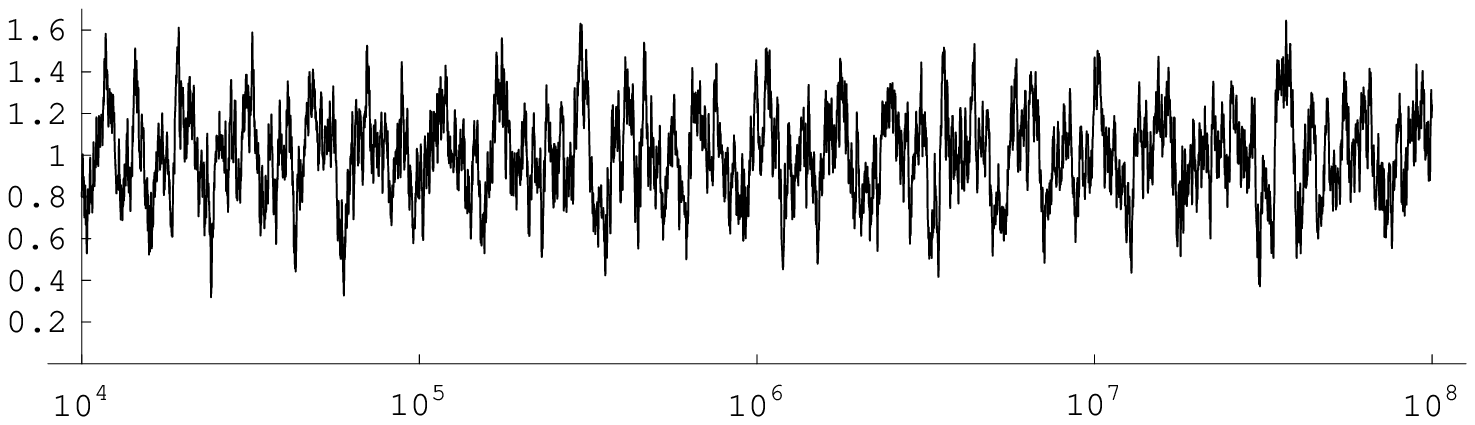} }
}
\botcaption{Figure 4} A graphs of $(\Li(x)-\pi(x)) \big/
(\frac12\Li(\sqrt{x}))$, followed by approximations using 10, 100, and 
1000
zeros of $\zeta(s)$.
\endcaption
\endinsert

    Extensive computations of this type  led Bays and Hudson
to conjecture that the first time that $\pi(x)$ surpasses $\Li(x)$
is at approximately $1.3982 \times 10^{316}$.

\subhead  The Mod 4 Race \endsubhead\smallskip

In 1959, Shanks suggested studying the Mod 4 Race by drawing a
histogram of the values of
$$
\frac{\#\{\text{primes}\ 4n+3 \le x\} - \#\{\text{primes}\ 4n+1
\le x\}}{\sqrt{x}/\ln x}. \tag{4}
$$
Figure 5 displays such a histogram for the one thousand sample
values $x={}$1000, 2000, 3000, \dots, $10^6$.

\midinsert {\centerline{
\epsfxsize=6truein\epsfbox{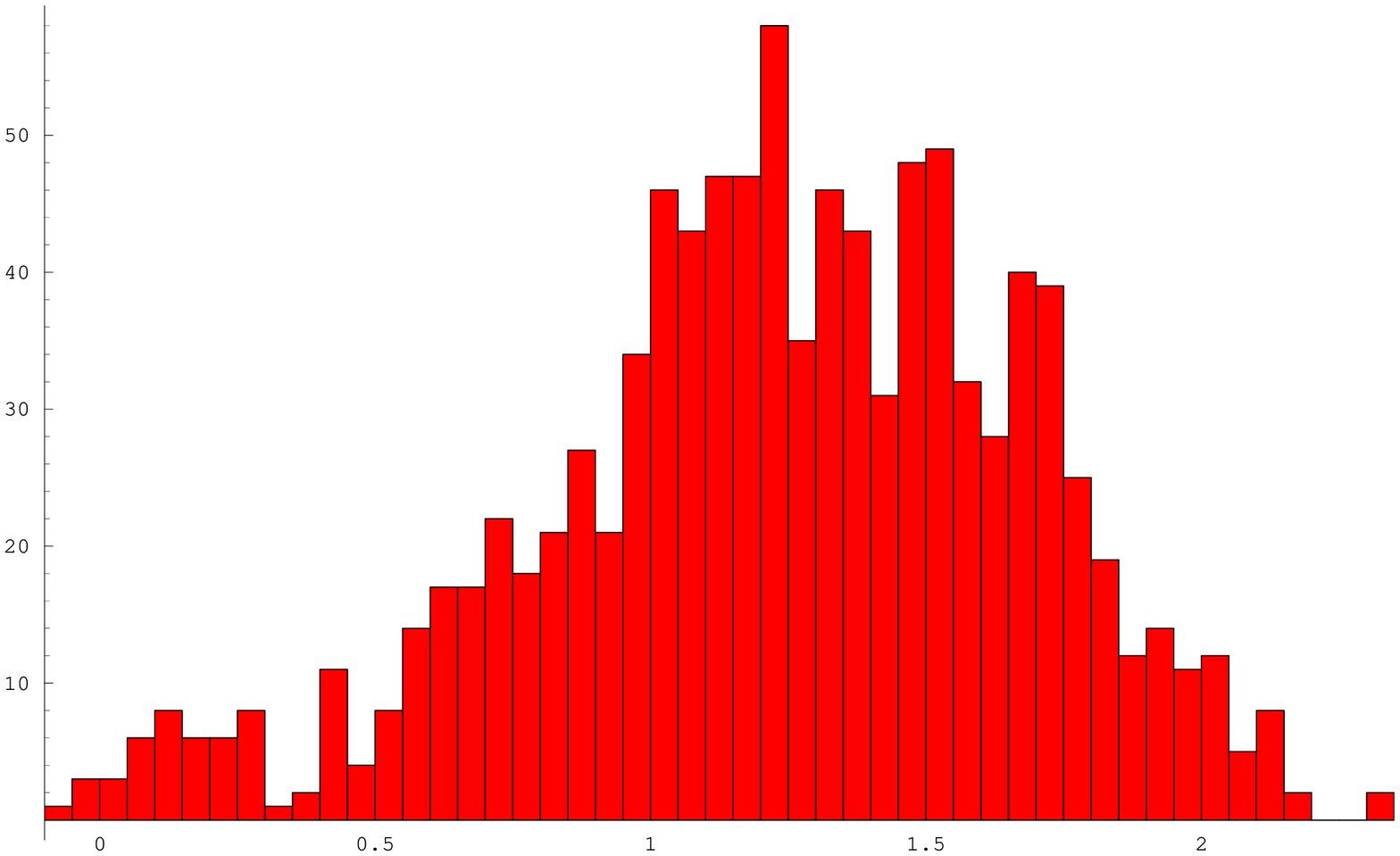} }} \botcaption{Figure
5} A histogram for the values of (4) at $x=1000k,\ 1\leq k\leq
1000$
\endcaption
\endinsert

\noindent This histogram is suggestive; one might guess that  if
we incorporate more and more sample values, then the histogram
will more and more resemble something like a bell-shaped curve
centered at~1.
    Littlewood's result (2) implies that the tail  does stretch out
horizontally to $\infty$ as more and more sample values are used,
since the ratio in equation (4) will be at least as large as
some constant multiple of $\ln\ln\ln x$ for infinitely many values of
$x$.\footnote{In
fact, Littlewood also proved the inequality with the two terms on
the left-hand side reversed, so that the histogram stretches out
to $-\infty$ as well.} The infinite extent of the histogram is
certainly not evident from the figure above, but this is not so
surprising, since

\smallskip
``$\ln \ln \ln x$ {\sl goes to infinity with great dignity}.''
\hfill --- Dan Shanks, 1959
\bigskip

In (3) we saw how the difference between $\pi(x)$ and $\Li(x)$ can
be approximated by a sum of waves whose frequencies and amplitudes
depend on the zeros of the Riemann zeta-function. To count primes
up to $x$ of the form $4n+1$, or of the form $4n+3$, or of any
form $qn+a$ with gcd$(a,q)=1$, there is a formula analogous to (3)
that depends on the zeros of {\it Dirichlet $L$-functions},
relatives of the Riemann zeta-function, which also have natural
though slightly more complicated definitions. For example, the
Dirichlet $L$-function associated to the race between primes of
the form $4n+3$ and primes of the form $4n+1$ is
$$
L(s) = \frac1{1^s} - \frac1{3^s} + \frac1{5^s} - \frac1{7^s} + \dots.
$$
Note that this (and all other Dirichlet $L$-functions) converges
for any $s=\sigma+it$ with $\sigma>0$. There is a lovely formula,
analogous to (3), for the number of primes of the form $qn+a$ up
to $x$. This formula holds if and only if all of the zeros of these
Dirichlet $L$-functions in the strip $0\leq \Re(s)\leq 1$ satisfy
$\Re(s)=1/2$. We call this assertion the {\sl Generalized Riemann
Hypothesis}. For example, if the Generalized Riemann Hypothesis is
true for the function $L(s)$ just defined, then we get the formula
$$
\frac{\#\{\text{primes}\ 4n+3 \le x\} - \#\{\text{primes}\ 4n+1
\le x\}} {\sqrt x/\ln x} \approx 1 + 2\dsize\sum \Sb \text{Real
numbers }  \gamma>0 \\ \text{such that } \frac 12 + i\gamma \\
\text{is a zero of }L(s)\endSb \ \ \frac{\sin(\gamma\ln x)}\gamma.
$$
Consequently, we can try to study the Mod 4 Race by computing the
right-hand side of this formula, truncating to involve only a
hundred or a thousand or a million zeros of $L(s)$. Here is a
graph showing how well such an approximation with 1000 zeros agrees
with the actual data:

\midinsert {
\centerline{ \epsfxsize=5truein\epsfbox{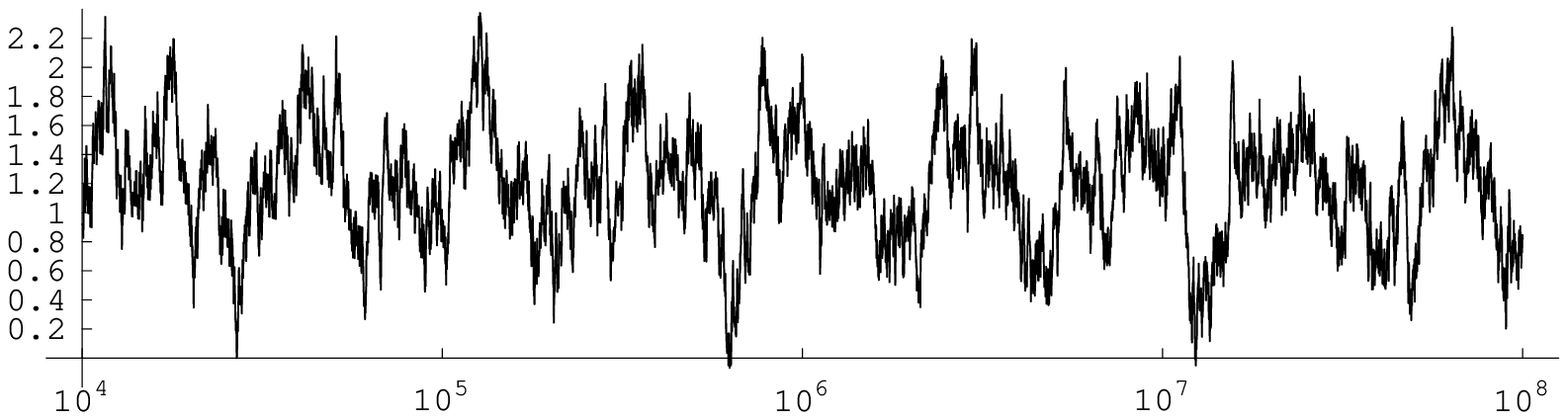} }
\vskip12pt
\centerline{ \epsfxsize=5truein\epsfbox{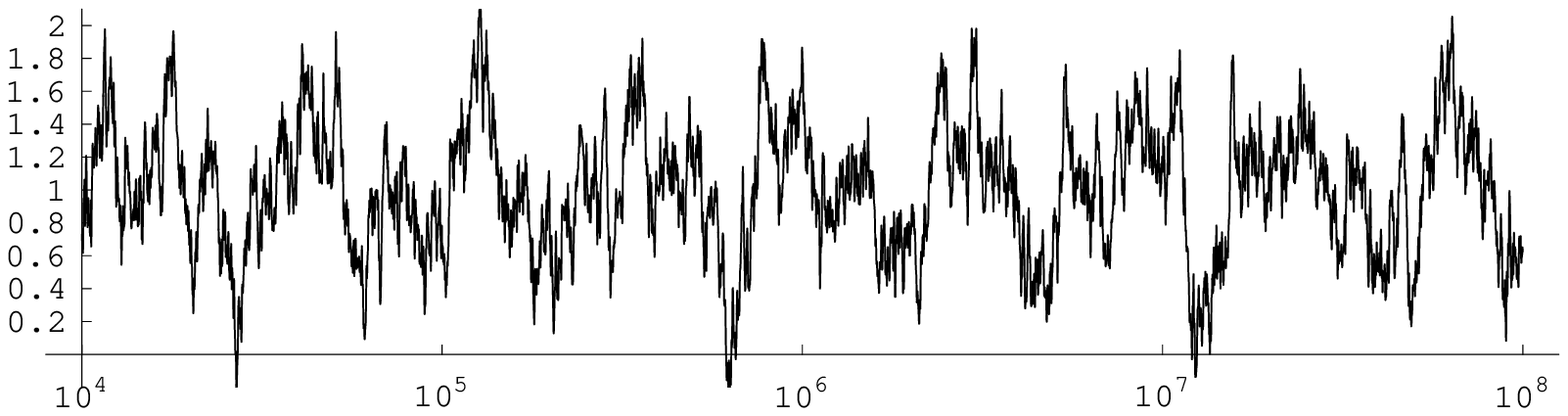} }
}
\botcaption{Figure 6} A graph of the function in (4), and an
approximation using 1000 zeros of $L(s)$.
\endcaption
\endinsert

The first three occasions that Team 1 takes the lead are clearly
visible in both graphs (above we gave the exact values at which
this happens). The convenient thing about the approximation is
that the approximation doesn't become hopelessly more difficult to
compute as
$x$ becomes large, unlike determining the precise count of primes.
In 1999 Bays and Hudson used such approximations to predict
further ``sign changes'' in (4) for $x$ up to $10^{1000}$. The
next two that they predicted were
$$
\align
& \approx 1.4898 \times 10^{12},\ \text{which was found at
1,488,478,427,089;} \\
& \approx 9.3190 \times 10^{12} . \\
\endalign
$$




\head Where do these biases come from? \endhead

We have observed several examples of prime number races and of
biases that certain arithmetic progressions seem to have over
others. Although we have seen that such bias can be calculated
through complicated formulae like (3), this explanation is not
easily understood, nor does it give us a feeling as to how we
might predict, without enormous calculation, which of two
progressions typically dominates the other. Let's summarize the
prime races we have seen so far and look for a pattern, hoping
that we can find some short cut to making such predictions. We
have seen that there seem to be, at least most of the time,
$$
\align
& \text{More primes of the form}\ 4n+3\ \text{than of the form}\ 4n+1;
\\
& \text{More primes of the form}\ 3n+2\ \text{than of the form}\ 3n+1;
\\
& \text{More primes of the form}\ 10n+3\ \text{and}\ 10n+7\ \text{than
of the form}\ 10n+1\ \text{and}\ 10n+9. \\
\endalign
$$
Do you see a pattern? There's probably not enough data here to
make a good guess. But let's check out one more four-way prime
number race, this one with the forms $8n+1, 8n+3, 8n+5, 8n+7$ as
the contestants:
$$
\matrix
\underline{x} & \underline{8n+1} & \underline{8n+3} & \underline{8n+5}
& \underline{8n+7} \\
1K & 37 & 44 & 43 & 43 \\
2K & 68 & 77 & 79 & 78 \\
5K & 161 & 168 & 168 & 171 \\
10K & 295 & 311 & 314 & 308 \\
20K & 556 & 571 & 569 & 565 \\
50K & 1,257 & 1,295 & 1,292 & 1,288 \\
100K & 2,384 & 2,409 & 2,399 & 2,399 \\
200K & 4,466 & 4,495 & 4,511 & 4,511 \\
500K & 10,334 & 10,418 & 10,397 & 10,388 \\
1\ \text{million} & 19,552 & 19,653 & 19,623 & 19,669 \\
\endmatrix
$$
\botcaption{Table 7} The number of primes of the form $8n+j$ for
$j=1,3,5$ and $7$ up to various $x$.
\endcaption
\medskip

\noindent  This is a strange race, since the only clear pattern
that emerges is that there seem to be

\smallskip
\centerline{More primes of the form $8n+3$, $8n+5$, and $8n+7$ than of
the form $8n+1$.}
\smallskip
\noindent (In fact, this is true for all $x$ between $23$ and
$10^6$.) But this ``strangeness'' should help students of number
theory guess at a pattern, because such students know that all odd
squares\footnote{Note that $(2m-1)^2=8\binom m2 +1\equiv 1 \pmod
8$ for all odd integers $2m-1$.} are of the form $8n+1$.

During our discussion about $\lcm[1,2,...,n]$ and its connection
with the prime number counting function $\pi(x)$, we saw that
Riemann's prediction $\pi(x)+\frac12\# \{ $primes $p$ with
$p^2\leq x\}$ is what is best approximated by $\int^x_2
\frac{dt}{\ln t}$, and so it is the squares of primes that
``account'' for the bias that $\Li(x)$ has over $\pi(x)$ for most
of the time.

These pieces of evidence\footnote{admittedly somewhat
circumstantial evidence} point to the squares of primes playing
an important role in such biases. Let's check out the arithmetic
progressions that the squares of primes belong to for the moduli
above.
$$
\align
\text{For any prime}\ p \ & \text{not\ dividing}\ 4,\ \text{we have}\
p^2\ \text{of\ the\ form}\ 4n+1; \\
\text{For any prime}\ p \ & \text{not\ dividing}\ 3,\ \text{we have}\
p^2\ \text{of\ the\ form}\ 3n+1; \\
\text{For any prime}\ p \ & \text{not\ dividing}\ 10, \ \text{we have}\
p^2\ \text{of\ the\ form}\ 10n+1\ \text{or}\ 9; \\
\text{For any prime}\ p \ & \text{not\ dividing}\ 8,\ \text{we have}\
p^2\
\text{of\ the\ form}\ 8n+1. \\
\endalign
$$
Exactly right every time! It does seem that ``typically'' $qn+a$
has fewer primes than $qn+b$ if $a$ is a square modulo $q$ while
$b$ is not.

\medskip \noindent {\bf What is the phenomenon that we are actually
observing?}\medskip

We have been observing some phenomena where there is a
``typical'' bias in favour of one arithmetic progression over
another, that is, one ``usual'' leader in the prime number race.
So far, though, we have not exhibited a precise description of
this bias nor defined a number that measures how strong the bias
is. In an earlier section, we stated Knapowski and Tur\'an's 1962
conjecture: if we pick a very large number $x$ ``at random'', then
``almost certainly'' there will be more primes $4n+3 \le x$ than
$4n+1 \le x$. (Obviously, such a conjecture can be generalized to
the other prime races we considered). However, the supporting
evidence   from our data was not entirely convincing. Indeed, by
studying the explicit formula (3), Kaczorowski (1993) and Sarnak
(1994), independently, showed that the Knapowski--Tur\'an
conjecture is false! In fact, the quantity
    $$
{\textstyle\frac 1X} \#\{x \le X \colon\ \text{there\ are\ more\ primes
of the form}\ 4n+3\ \text{up to $x$ than of the form}\ 4n+1\}
$$
does not tend to any limit as $X \rightarrow \infty$, but  instead
fluctuates. This is not a very satisfying state of affairs, though
it is feasible that we can't really say anything more: it could be
that the phenomena we have observed are only true for ``small
numbers'',\footnote{``small'' on an appropriate scale!} and that
the race looks truly ``unbiased'' when we go out far enough.

However, this turns out to not be the case. In 1994, Rubinstein
and Sarnak made an inspired observation. After determining that
one progression does not dominate the other for any fixed
percentage of the time,
so that that precise line of
thinking seemed closed, Rubinstein and Sarnak said to themselves
that perhaps the {\it obvious\/} way to count prime number races
is not necessarily the {\it correct\/} way to count prime number
races, at least in this context. They observed  that if you count
things a little differently (in technical terms, ``if you use a
different measure''), then you get a much more satisfactory
answer. Moreover, you get a whole panorama of results that explain
all the observed phenomena in a way that feels very natural.

Sometimes in mathematics it is the simplest remark that is the
key to unlocking a mystery, not the hard technical slog, and that
was certainly the case for this question. The key idea of
Rubinstein and Sarnak is to not count 1 for each integer $x\leq X$
for which there are more primes of the form $4n+3$ up to $x$ than
of the form $4n+1$, but rather to count $1/x$.  Of course, the
total sum $\sum_{x\leq X} 1/x$ is not $X$, but rather
approximately $\ln X$, so we need to scale with this instead.
However, when we do so, we obtain a remarkable result:

\medskip
\noindent {\bf Rubinstein and Sarnak} (1994):
$$
\frac{1}{\ln X}\ \sum \Sb x \le X \\ \text{There are more\ primes}\
4n+3 \\
\text{than}\ 4n+1\ \text{up\ to}\ x \\ \endSb\ \frac{1}{x}\
\longrightarrow \
.9959\ldots\ \text{as}\ X \to \infty .
$$
In other words, with the appropriate method of measuring, Tch\'ebychev
was correct
99.59\% of the time!

Moreover their idea can be  applied to all of the other  prime
number races that we have been discussing. For example, in the Mod
3 Race, we have
$$
\frac{1}{\ln X}\ \sum \Sb x \le X \\ \text{There are more\ primes}\
3n+2 \\
\text{than}\ 3n+1\ \text{up\ to}\ x \\ \endSb\ \frac{1}{x}\
\longrightarrow \
.9990\ldots\ \text{as}\ X \to \infty,
$$
so that Team 1 only has about one chance in a thousand of being
ahead of Team 3 at any given time, when we measure ``being ahead''
this way.

And what about the race between $\pi(x)$ and $\Li(x)$? Remember that we
don't expect to find a
counterexample until we get out to the huge value $x\approx 1.3982
\times 10^{316}$. In this case,
Rubinstein and Sarnak showed that
$$
\frac{1}{\ln X} \sum \Sb x \le X \\ \pi(x) < \text{Li}(x)  \endSb\
\frac{1}{x} \
\longrightarrow \ .99999973\ldots\ \text{as}\ X \to \infty .
$$
That is, with the appropriate method of measuring, $\pi(x)
<\Li(x)$  about $99.999973\%$ of the time! No wonder that the
first point at which $\pi(x) > \Li(x)$ is so far out.

This successful way of measuring things is called the {\sl
logarithmic measure} (because of the ``normalizing factor'' $\ln
X$). It appears in many contexts in mathematics, and was even used
in a related context by Wintner in 1941, but never quite in such
an appealing and transparent way as by Rubinstein and Sarnak.

Their wonderful results are proved under highly plausible
assumptions. They assume first that the Generalized Riemann
Hypothesis holds (in other words, marvellous formulae analogous to
(3) can be used to count primes of the form $qn+a$ up to $x$), and
second that the $\gamma$'s that arise in these formulae---the
imaginary parts of the zeros of the associated Dirichlet
$L$-functions---are actually linearly independent over the
rational numbers. It would be shocking if either of these
assumptions is false, but it does seem unlikely that either one
will be proved in the near future.

\medskip
\noindent{\bf  How do you prove such results?}
    \smallskip

The proofs of the results mentioned above, and those discussed in
the next sections, are too technical to describe in full detail
here. However, they do all depend on careful examination of the
formula (3) and its analogues for the Mod $q$ Races. The
assumption of the Generalized Riemann Hypothesis ensures, in all
the cases under consideration, that  a formula of the type (3)
exists in the first place. To analyze its value, we need to
consider the sum on the right side of (3) and, in particular, the
values of the various terms and how they interact: Each term
$-2\sin (\gamma \ln x)/\gamma$ is a sine-wave which undulates
regularly\footnote{``Regularly'', that is, if we consider $\ln x$
to be the variable rather than $x$. In fact, this is really the
reason that the ``logarithmic measure'' is most appropriate for
our prime number counts.} between $-2/\gamma$ and $2/\gamma$. As
we have seen, when studying how such waves can combine to give
approximations to a straight line, if the various values of
$\gamma$ can be chosen so that somehow the individual undulations
can be synchronized, then the sum can accumulate into surprising
shapes. This leads to the second assumption: If we assume that the
waves cannot combine in some extraordinary way, then the sum
should remain small and unsurprising. In other words, we would
expect that a fair portion of the terms would be positive and a
fair portion will be negative, so that there will be significant
cancellation. This is what is achieved by assuming that the
$\gamma$'s are linearly independent over the rational numbers.
Under this second assumption, the values of the $\sin (\gamma \ln
x)$ terms cannot be well synchronized for very many values of $x$,
and so we are able to prove results. In fact, we find that the
smallest few $\gamma$'s have the largest influence on the value of
the sum (which is not surprising, since the term corresponding to
$\gamma$ has a factor of $\gamma$ in the denominator).

Calculating precise numerical values for how often one prime
number count is ahead of another is even more delicate. It is no
longer enough to argue that the sum is ``large'' or ``small''; we
need to know exactly how often the sum in (3) is greater than
$-\frac12$ or less than $-\frac12$, since that value is where the
right-hand side switches between positive and negative. As it
turns out, it is possible to calculate exactly how often this
occurs if we use (traditional) Fourier analysis: We rewrite the
sum of the waves in (3) in which the frequencies depend on the
zeros of the Riemann zeta-function, as a sum of ``Fourier waves''
where the frequencies depend on the numbers $2\pi n$. Although
this is too cumbersome to give in a closed form, one can, with
clever computational techniques, approximate its numerical value
    with great accuracy.

\subhead Squares and nonsquares:\ Littlewood's method
\endsubhead\smallskip

On several occasions above we encountered consequences of
Littlewood's 1914 paper\footnote{Which actually only dealt with
the $\pi(x)$ vs.~ $\Li(x)$ race.}. The method is robust enough to
prove some fairly general results but not typically about one
progression racing another\footnote{Though that was all we saw in
the cases above where $q$ is small.}. Indeed, for any odd prime $q$, let
Team S be the set of primes that are congruent to a square mod
$q$, and Team N be the rest of the primes. Using Littlewood's
method one can show that each such team takes the lead infinitely
often\footnote{In other words, there are arbitrarily large $x$ and
$y$, such that there are more primes in $S$ than in $N$ up to $x$;
and there are more primes in $N$ than in $S$ up to $y$.}, and the
lead up to $x$ can be as large as $c\sqrt{x} \ln\ln\ln x/\ln x$
for some constant $c>0$ that depends only on $q$. One  example of
this is in the $q=5$ race where $S$ contains the primes with last
digits 1 and 9 (as well as 5), so this is the Mod 10 Race
again---recall that we observed in the
limited data above that Team N typically held the lead.

There are other races that can be analyzed by Littlewood's method,
always involving partitioning the primes into two seemingly equal
classes. In fact there is at least one such race for every modulus
$q$, and often more than one. When $q$ is not a prime, however,
the descriptions of the appropriate classes become more
complicated.




    \subhead More results from Rubinstein and Sarnak
(1994)\endsubhead\smallskip

In a prime number race between two arithmetic progressions modulo
$q$, when do we see a bias and when not? Is each arithmetic
progression ``in the lead'' exactly 50\% of the time (in the
logarithmic measure) or not? More importantly perhaps, can we
decide this before doing a difficult calculation? Above we saw
that there seem to ``usually'' be more primes up to $x$ of the
form $qn+b$ than of the form $qn+a$ if $a$ is a square modulo $q$
and $b$ is not. Indeed, under the two assumptions mentioned above,
Rubinstein and Sarnak proved that this is true: the logarithmic
measure of the set of $x$ for which there are more primes of the
form $qn+b$ up to $x$ than of the form $qn+a$  is strictly greater
than $\frac12$, although always less than~1. In other words, any
nonsquare is ahead of any square more than half the time though
not 100\% of the time.

    We can ask the same question when $a$ and $b$ are either
both squares modulo $q$ or both nonsquares modulo~$q$. In this
case, under the same assumptions, Rubinstein and Sarnak proved
that $\#\{\text{primes}\ qn+a \le x\}> \#\{\text{primes}\ qn+b \le
x\}$ exactly half the time. In fact they prove rather more than
this. To describe their result, we need to define certain error
terms related to these prime number counts and describe what we
mean by their ``limiting distributions''. As we noted above, the
values of the prime counting function $\#\{\text{primes}\ qn+a \le
x\}$ are all roughly equal as we vary over the integers $a$ up to
$q$ that have no factor in common with $q$. In fact, we saw that
the ratio of any two of them tended to 1 as $x\to \infty$. This
implies that
$$
\frac{ \#\{\text{primes}\ qn+a \le x\} }{ \pi(x)/\phi(q)}
\to 1 \qquad \text{as} \ x\to \infty,
$$
where $\pi(x)=\#\{\text{primes} \le x\}$ as before and $\phi(q)$ is the
number of positive integers $a$ up to $q$ for which $\gcd(a,q)=1$.
We have seen that it is natural to look at the difference in these
quantities divided by $\sqrt{x}/\ln x$, so let us define
$$
\text{Error}(x;q,a) = \frac {  \#\{\text{primes}\ qn+a \le x\} -
\pi(x)/\phi(q) } { \sqrt{x}/\ln x }.
$$

Rubinstein and Sarnak suggested that to  properly study prime
races between arithmetic progressions $a\pmod q$ and $b\pmod q$,
one should look at the distribution of values of the ordered pair
$$
\big( \text{Error}(x;q,a), \text{Error}(x;q,b) \big)
$$
as $x$ varies, this distribution again being  defined with respect
to the logarithmic measure. More concretely, given real numbers
$\alpha<\beta$ and $\alpha'<\beta'$, one might ask how often
$\alpha\leq \text{Error}(x;q,a) \leq \beta$ and $\alpha'\leq
\text{Error}(x;q,b) \leq \beta'$. This frequency is defined by the
limit of integrals
$$
\lim_{Y\to \infty} \ \frac 1{\ln Y} \ \int\Sb 0\leq y\leq Y, \\
\text{Error}(y;q,a) \in [\alpha,\beta] \\
\text{and}\ \text{Error}(y;q,b) \in [\alpha',\beta']\endSb
\ \frac 1y \,dy
$$
(a limit that they prove exists). Rubinstein and  Sarnak proved
that if $a$ and $b$ are both squares or both nonsquares modulo
$q$, then this distribution is symmetric, that is,
$$\alpha\leq \text{Error}(x;q,a) \leq \beta \quad\text{and}\quad
\alpha'\leq \text{Error}(x;q,b) \leq \beta'
$$
as often as
$$
\alpha'\leq \text{Error}(x;q,a) \leq \beta' \quad\text{and}\quad \alpha
\leq \text{Error}(x;q,b) \leq \beta.
$$
In other words, there is no sign of any bias, of any form,  at
all: the arithmetic progressions $a\pmod q$ and $b\pmod q$ are
interchangeable in the limit.
\smallskip

This was all studied in much more generality. For example,  one
can ask about the twenty-four possible orderings of the four prime
number counts
$$
\#\{\text{primes}\ 8n+1 \le x\}, \#\{\text{primes}\ 8n+3 \le x\},
\#\{\text{primes}\ 8n+5 \le x\}, \#\{\text{primes}\ 8n+7 \le x\}.
$$
Rubinstein and Sarnak showed that each ordering occurs for
infinitely many values of $x$; in fact, each ordering occurs a
positive proportion of the time (in the logarithmic measure). This can
be
generalized to all moduli $q$ and to as many distinct arithmetic
progressions $a_1,\dots a_r \pmod q$  as we like. That is, one can
study the distribution of
$$
(\text{Error}(x;q,a_1), \text{Error}(x;q,a_2), \dots,
\text{Error}(x;q,a_r)).
$$
Assuming only the Generalized Riemann Hypothesis, they prove that
the distribution function\footnote{Defined by a certain limit of
integrals, analogous to the limit a few lines earlier.} for this
vector of error terms exists. Assuming also the linear
independence of the $\gamma$'s over the rational numbers, they
prove that  for any distinct $a_1,\dots a_r \pmod q$, none of
which have any factors in common with $q$, the ordering
$$
\#\{\text{primes}\ qn+a_1 \le x\} < \#\{\text{primes}\ qn+a_2 \le x\} <
\dots < \#\{\text{primes}\ qn+a_r \le x\}
$$
occurs for a positive proportion (under the logarithmic  measure)
of values $x$. They also proved that, for any fixed $r$, these
proportions get increasingly close to $1/r!$ as $q$ gets larger.
It seems extremely unlikely that this proportion is usually
exactly $1/r!$, but we cannot prove that it is or isn't except in
the following special situation.

Above we saw that the distribution function is symmetric  when we
have a two-progression race and the progressions are either both
squares or both nonsquares modulo $q$. Surprisingly, Rubinstein
and Sarnak showed that there is only one other situation in which
the distribution function is symmetric no matter how we swap the
variables with one another, namely the race between three
arithmetic progressions of the form
$$
a \pmod q,\ a\omega \pmod q, \ \text{and} \  a\omega^2 \pmod q
$$
where $\omega^3\equiv 1 \pmod q$ but $\omega \not\equiv 1 \pmod q$.

However, Rubinstein and Sarnak's result that the distribution
function is not usually symmetric still leaves open the
possibility that each ordering in a race occurs with the same
frequency\footnote{Since (for example) a function can be positive
half the time without being symmetric.}.

Earlier we saw Shanks' histogram of values of
    $\text{Error}(x;4,3)-\text{Error}(x;4,1)$ for various values of $x$,
and guessed
    that ``the histogram will more and more resemble something like a
bell-shaped curve centered at~1'' as we take more data points. One
perhaps surprising consequence of the work of Rubinstein and
Sarnak is that this most obvious guess is not the correct one:
Although there is a distribution function, it will not be as
simple and elegant as the classical ``normal distribution'' bell
curve.

    \head What research is happening
right now? \endhead


After Littlewood's great (and, at the time, quite surprising) results
in 1914, there was a lull in research in
    ``comparative prime number theory'' until the fifties and sixties. At
that time, Littlewood's ideas
were extended in several theoretical directions that were arguably
suggested directly by Littlewood's work, and by substantial
calculations of Shanks, Hudson, and others. By the nineties, it
appeared that most of what could be done in this subject had been
done (and that much would never be done), although Kaczorowski was
still proving some new results at this time.

In 1993, Giuliana Davidoff taught a senior-level undergraduate
course at Mount Holyoke College in analytic number theory. Finding
so many interested students, Professor Davidoff decided that she
could run a fun ``{\sl Research Experience for Undergraduates}''
(REU) that summer investigating prime number races\footnote{Under
the more formal title ``{\sl Topics in comparative number
theory}''.}. Along with students Caroline Osowski,
Jennifer vanden Eynden, Yi Wang, and Nancy Wrinkle, she did some
important
calculations and proved several results that will be given in the
first appendix to this article.

By chance, Davidoff met Sarnak on summer vacation. She described
her REU project and how they were developing the computational
approach of Stark to such problems. Sarnak writes\footnote{In an
email to the first-named author.}: ``I was not familiar with this
Chebyshev bias feature and became fascinated. In particular I was
\dots very interested to put a definite number as to the
probability of $\pi(x)$ beating $\Li(x)$.''

Sarnak continues: ``When I got back to Princeton I chatted with
Fernando\footnote{Fernando Rodriguez-Villegas, now of the
University of Texas at Austin, but then a Princeton postdoctoral
fellow who had lots of computational experience.} about this
\dots and began to work on this \dots It was clear that many
zeroes would have to be computed''. Sarnak had previously
discussed other questions on the distribution of prime numbers
with a bright undergraduate, Mike Rubinstein. ``Mike, who was
looking for a senior thesis topic, got very involved in this and
after a while he and I worked on the problem as collaborators and
this led to our paper.'' This was quite an extraordinary senior
thesis, as it became one of the most influential papers in recent
analytic number theory. Subsequently Rubinstein went on to his
Ph.D.\footnote{Michael Rubinstein is now an assistant professor at
the University of Waterloo.}~and has become one of the world's
leading researchers in the computation of different types and
aspects of zeta-functions.

Soon thereafter, the second
author read Rubinstein and Sarnak's paper and
became   interested in determining the ``probabilities'' for some
of the three-way prime number races that Rubinstein and Sarnak's
work did not address---for example, the race between the three
contestants
$$
\#\{\text{primes } 8n+3\le x\},\, \#\{\text{primes } 8n+5\le x\}\
\text{and } \#\{\text{primes } 8n+7\le x\};
$$
note that none of these arithmetic progressions contain any
squares.

With colleague Andrey Feuerverger at the University of
Toronto\footnote{Where he was a postdoctoral fellow at the time before
moving to the University of British Columbia.}, he created a technique
to
determine the frequencies of the various orderings of contestants
in a prime number race with more than two contestants.
Unexpectedly they found that the six orderings of  the three
contestants in the race above do not necessarily each happen
one sixth of the time. In fact
$$
\#\{\text{primes}\ 8n+3 \le x\} > \#\{\text{primes}\ 8n+5 \le x\} >
\#\{\text{primes}\ 8n+7 \le x\}
$$
holds for approximately 19.2801\% of integers $x$ (in the
logarithmic measure), while
$$
\#\{\text{primes}\ 8n+5 \le x\} > \#\{\text{primes}\ 8n+3 \le x\} >
\#\{\text{primes}\ 8n+7 \le x\}
$$
holds for approximately 14.0772\% of integers $x$ (under
appropriate assumptions).

This is strange, for in the  race between primes of the form
$8n+3$ and primes of the form $8n+5$ both contestants are in the
lead half the time, yet if we {\it only\/} look at values of $x$
for which both of these prime number counts are ahead of the count
for primes of the form $8n+7$, then the $8n+3$ primes are more
likely to be ahead!

\subhead More esoteric races \endsubhead

There are results known that give asymptotics for the number of
primes having certain properties. For example, when is 2 a cube
modulo the prime $p$? It is easy to see that this is the case for all
primes of the form
$3n+2$, and so we focus on the other primes. It is known that 2 is a
cube for one-third of the
primes of the form $3n+1$, asymptotically.  So we ask whether,
typically,
slightly more or slightly less than one-third of the primes $p=3n+1$
have the property that $2$ is a cube modulo~$p$.

In his recent doctoral thesis at the University
of British Columbia, Nathan Ng\footnote{Nathan Ng is now an assistant
professor at the University of Ottawa.} noted that this question can be
answered using similar techniques to Rubinstein
and Sarnak, since the count of these primes depends in an
analogous way on the zeros of yet another type of $L$-function.
Indeed Ng observed that such results  can be generalized to many
cases in which we know the asymptotics for a set of primes that
can be described by Galois theory\footnote{The $L$-functions in
question are Artin $L$-functions, and the asymptotics for
conjugacy classes under the appropriate Frobenius maps are
obtained by Cebotarev's density theorem. Analogous results are
proved under the assumption of holomorphy for these $L$-functions
as well as the Riemann Hypothesis and linear independence over the
rationals of the zeros of these $L$-functions.}.

Ng gave the following nice example of his work.  Power series like
$$
q \prod_{n=1}^{\infty}(1-q^{n})(1-q^{23n}) = \sum_{n=1}^{\infty}
a_{n} q^{n}
$$
appear frequently in arithmetic geometry; this is an example of a
{\sl modular form}, and has all sorts of seemingly miraculous
properties (see Serre's book [Se]). One such miracle is that for
every prime $p$, the value of $a_p$ is $2,0,$ or $-1$ and that the
proportions of each are around $1/6, 1/2$ and $1/3$, respectively.
Under the appropriate assumptions Ng showed, in the logarithmic
measure, that

$2\#\{p \le x, a_{p}=0\}>6\#\{p \le x, a_{p}=2\}$ for $\approx 98.30\%$
of the values of $x$,

$2\#\{p \le x, a_{p}=0\}>3\#\{p \le x, a_{p}=-1\}$ for $\approx
72.46\%$ of the values of $x$, and

$3\#\{p \le x, a_{p}=-1\}>6\#\{p \le x, a_{p}=2\}$ for $\approx
95.70\%$ of the values of $x$.

\subhead Prime pairs \endsubhead\smallskip

    The first author's
involvement in this topic came from an invitation to speak at an
MAA meeting in Montgomery, Alabama in 2001, which gave him an
excuse to read up on this fascinating subject. A subsequent
discussion with students led to the idea of creating a 2001-2002
VIGRE\footnote{The latest acronym to come out of Washington, a new
type of grant which encourages ``Vertical Integration'' of
research, from bright undergraduates through senior faculty.}
research group at the University of Georgia to investigate
analogous questions about twin prime races:

For any even number $d$, we believe that there are infinitely many
integers $n$ for which both $n$ and $n+d$ are prime. When $d=2$,
this is the famous {\sl Twin Prime Conjecture}. At present, we
cannot prove that there are infinitely many ``prime pairs'' $n,
n+d$ for any value of $d$, whatsoever. Nonetheless, we can predict
how many there should be up to $x$. Indeed, we believe that
$$
\frac{ \#\{ \text{prime pairs }n,n+2^j \le x \} }{ \#\{ \text{prime
pairs }n,n+2^k \le x \} } \to 1 \quad \text{as } x\to\infty
$$
for any positive integers $j$ and $k$. Comparing the first four
powers of~2, we obtain the following data:
$$
\matrix
\underline{n \le x} & \underline{n,n+2} & \underline{n,n+4} &
\underline{n,n+8} & \underline{n,n+16} \\
100 & 8 & 9 & 9 & 9 \\
200 & 15 & 14 & 14 & 13 \\
500 & 24 & 27 & 24 & 24 \\
1000 & 35 & 41 & 38 & 39 \\
2000 & 61 & 65 & 63 & 60 \\
5000 & 143 & 141 & 141 & 135 \\
10000 & 205 & 203 & 208 & 200 \\
20K & 342 & 344 & 353 & 331 \\
50K & 705 & 693 & 722 & 707 \\
100K & 1224 & 1216 & 1260 & 1233 \\
200K & 2160 & 2136 & 2194 & 2138 \\
500K & 4565 & 4559 & 4641 & 4631 \\
1\ \text{million} & 8169 & 8144 & 8242 & 8210 \\
\endmatrix
$$
\botcaption{Table 8} The number of prime pairs $n,\ n+2^k$ up to
$x$,   for $k=1,2,3,4$ and various values of $x$.
\endcaption
\medskip

It appears from this initial data that there is a bias in favor of
prime pairs $n,n+8$, at least once $x$ is ten thousand or more.
However, there is no explicit formula like (3) known that would
lead to a  fruitful attack on this question. In the second
appendix to this paper, we present the work of a team of students
at the University of Georgia that sought further data and tried to
make predictions.

\smallskip

There are several other natural questions pertaining to prime
number races that have been the subject of research over the last
few years.

\subhead It can take a long time to take the lead
\endsubhead\smallskip

Although Team 1 does occasionally get ahead in the Mod 4 Race, it
takes quite a while for this to happen. Indeed $x$ is larger than
25,000 the first time Team 1 takes the lead and then only
momentarily; the first time Team 1 sustains the lead over
a longer interval $x$ is larger than half a million. Similarly, in
the Mod 3 Race,  $x$ is larger than half a trillion before Team 1
takes the lead. In the $\pi(x)$ versus $\Li(x)$ race, we still do
not know exactly when $\pi(x)$ first takes the lead though the
answer is definitely gigantic.

Another situation is where one Team is ganged up on by the other
Teams. For instance we saw that primes of the form $8n+1$ lag far
behind primes of the forms $8n+3,\ 8n+5$ or $8n+7$ (because all
the odd squares, including squares of primes, are of the form
$8n+1$, a heavy burden on Team 1). In fact, the first time that
Team 1 even gets out of last place is beyond half a billion:
$x=588,067,889$ is the first time that there are more primes of
the form $8n+1$ up to $x$ than there are primes of the form
$8n+5$; then $x\approx 1.9282\times 10^{14}$ is the first time
that there are more primes of the form $8n+1$ up to $x$ than there
are primes of the form $8n+7$.  Kaczorowski and, independently,
Rubinstein and Sarnak showed that Team 1 runs in first place in
this four-way race for a positive proportion of values of $x$;
however, the first time this occurs is beyond $10^{28}$.

So why does it take so long for some competitors to get their turn
in the lead? To understand this we need to examine the most
``important'' terms in formulas like (3), namely the ``1'' term
and the first few summands, those for which $\gamma$ is small.
Looking closely at the formula (3), these first few summands,
which correspond to the nontrivial zeros of the Riemann zeta
function that are closest to the real axis, have $\gamma$
approximately equal to 14.13, then 21.02, 25.01, 30.42, \dots. In
particular, these summands are all $<1/7$ in absolute value, which
is small relative to the ``1'' term. Even if we could find a value
of $x$ for which many of the initial terms $\sin(\gamma\ln x)$ are
close to $-1$, it would take at least the first twenty-one terms
of the sum to be negative enough to compensate for the initial
``1'' term. Thus these values of $x$ must be special in that many
of the $\sin(\gamma\ln x)$ terms must be close to $-1$
simultaneously.

If we look at the Mod $q$ Race as $q$ gets larger, it turns out
that the relevant numbers $\gamma$ become smaller (that is, the
zeros of the appropriate Dirichlet $L$-functions lie closer to the
real axis), and therefore  changes in the lead happen sooner and
more frequently.

However, tiny values of $\gamma$ come with their own problems. For
example, when $q=163$  the relevant sum has an especially small
$\gamma$, namely $\gamma\approx 0.2029$. So this one summand has a
large impact
on the final answer since the denominator is so small. The first
time $\sin (\gamma \ln x)$ is close to $-1$ is when $\gamma \ln x$
is close to $3\pi/2$; which corresponds to a value of $x$ around
12 billion! And if that $x$ doesn't work out because of the other
terms then the next such value is when $\gamma \ln x$ is close to
$7\pi/2$, which corresponds to $x\approx3.43\times10^{23}$. No
wonder it can take a while to see the predicted effects!

\subhead Sharing the lead \endsubhead\smallskip

In a race between two contestants who keep taking the lead from
one another, there will be plenty of moments when the two
contestants are tied: Since there are arbitrarily large values of
$x$ for which $\#\{ \text{primes}\ qn+a\leq x\} > \#\{
\text{primes}\ qn+b\leq x\}$, and also arbitrarily large values
for which $\#\{ \text{primes}\ qn+a\leq x\} < \#\{ \text{primes}\
qn+b\leq x\}$, and since these counting functions take only
integer values, there must be infinitely many integers $x$ for
which $\#\{ \text{primes}\ qn+a\leq x\} = \#\{ \text{primes}\
qn+b\leq x\}$.

What about races with three or more contestants? Even if each of
the contestants leads at some point, there seems to be no particular
reason for all three of them to be tied at any point. So we ask whether
there
exist infinitely many integers $x$ for which
$$
\#\{ \text{primes}\ qn+a\leq x\} = \#\{ \text{primes}\ qn+b\leq x\} =
\#\{ \text{primes}\ qn+c\leq x\} = \cdots.
$$
Feuerverger and Martin conjecture that there {\it are\/}
infinitely many such ties in a three-way Mod $q$ Race, but {\it not\/}
in
a Mod $q$ Race with four or more teams. Their compelling argument
runs as follows: Consider the $(k-1)$-dimensional vector whose
$i$th entry is the difference
$$
\#\{ \text{primes}\ qn+a_i\leq x\} - \#\{ \text{primes}\
qn+a_{i+1}\leq x\} \tag{5}
$$
    for each $1\leq i\leq k-1$. Notice that the $k$ counting
functions $\#\{ \text{primes}\ qn+a_{i}\leq x\}$ are all tied with
one another precisely if this $(k-1)$-dimensional vector equals
$(0,0,\dots,0)$. As we let $x$ increase, this vector changes each
time $x$ equals a prime number in one of the arithmetic
progressions we are counting, changed by adding one of the vectors
$$
(1,0,\dots,0),\, (-1,1,0,\dots,0),\, \cdots,\,
(0,\dots,0,-1,1,0,\dots,0),\, \cdots,\, (0,\dots,0,-1), \tag{6}
$$
depending on whether the prime is of the form $qn+a_1$ or $qn+a_2$
or ... $qn+a_k$, respectively. The prime number theorem for
arithmetic progressions tells us that each of these vectors is
roughly equally likely to occur.

  From this, Feuerverger and Martin suggest that the progress of the
$(k-1)$-dimensional vector in (5) can be modelled by a ``random
walk'' on the $(k-1)$-dimensional lattice generated by the vectors
listed in (6). Now it is known that with probability 1, a random
walk on a one- or two-dimensional lattice will return to the
origin (indeed, will visit every lattice point) infinitely often,
but a random walk on a lattice of dimension three or greater will
return to the origin only finitely often. This is the source of
Feuerverger and Martin's conjecture, since the lattices of
dimension one and two in this model correspond to prime number
races with two and three contestants, respectively.

\subhead Certain special symmetries \endsubhead\smallskip

Assuming the Generalized Riemann Hypothesis, Feuerverger and
Martin proved  that certain configurations
$$
\#\{\text{primes}\ qn+a_1 \le x\} < \#\{\text{primes}\ qn+a_2 \le x\} <
\dots < \#\{\text{primes}\ qn+a_r \le x\}
$$
occur just as often\footnote{With respect to the logarithmic
measure of the set of such values $x$.} as certain other
configurations
$$
\#\{\text{primes}\ qn+b_1 \le x\} < \#\{\text{primes}\ qn+b_2 \le
x\} < \dots < \#\{\text{primes}\ qn+b_r \le x\}
$$
in the following situations:

\smallskip\noindent $\bullet$ The $a$'s and $b$'s are  inverses of each
other modulo $q$; that is, $a_ib_i\equiv 1 \pmod q$ for all~$i$;

\smallskip\noindent $\bullet$ The list of $b$'s is just the reversal of
the list of $a$'s, namely $b_i=a_{r+1-i}$ for all $i$,
and are either all squares modulo $q$ or all nonsquares modulo
$q$;

\smallskip\noindent $\bullet$ There exists an integer $m$ such that
$a_i \equiv mb_i \pmod q$ for each $i$, and in addition,
one of the following holds:

\smallskip$\bullet$ The $a$'s are all squares modulo $q$; or

$\bullet$ For each $i$, the two numbers $a_i$ and $b_i$ are either both
squares modulo $q$ or both nonsquares modulo $q$; or

$\bullet$ For each $i$, exactly one of $a_i$ and $b_i$ is a square
modulo $q$.
\smallskip

Probably the orderings appear with different frequencies if they
are not related by some special symmetry.

\medskip
\noindent {\bf  And what if the Riemann Hypothesis is
false?}
\smallskip

Certain famous results in number theory have been proved in two
parts: first under the assumption that the Riemann Hypothesis is
true, and second under the assumption that the Riemann Hypothesis
is false. The groundbreaking work of Sarnak and Rubinstein
followed from two hypotheses, the first of which was the
Generalized Riemann Hypothesis. One might ask whether we can
obtain similar results, perhaps with a rather different proof,
under the assumption that the Generalized Riemann Hypothesis is
false.

If the Generalized Riemann Hypothesis is false, it turns out to be
easier to prove that $\#\{ \text{primes}\ qn+a\leq x\} > \#\{
\text{primes}\ qn+b\leq x\}$ for a positive proportion\footnote{In
the logarithmic measure.} of $x$, since the ``1'' term in the
formula analogous to (3) becomes irrelevant, and this was  what
caused the bias when we assumed the Generalized Riemann
Hypothesis. However, whether one can prove that this happens 50\%
of the time is still in doubt. If one could prove this under the
assumption that the Generalized Riemann Hypothesis is false, then
we  would have  an unconditional proof that the Mod $q$ Race
between two squares or two nonsquares is always an even race (by
combining such a proof with the work of Rubinstein and Sarnak).

The behavior of prime number races might be genuinely different if
the Generalized Riemann Hypothesis does not hold. In 2001, Ford
and Konyagin proved that if the Generalized Riemann Hypothesis is
false---false in an absurdly convenient, yet feasible, way---then
there are some orderings of the prime number counting functions
that never occur.

To describe one of their constructions, we need to discuss a
generalization of the Riemann zeta-function. Define $\chi$ to be
that function for which
$$
\chi(5m)=0,\ \chi(5m\pm 1)=\pm 1,\, \text{and } \chi(5m\pm 2)=\pm
i \ \ \text{for all integers}\ m;
$$
and let $L(s,\chi)=\sum_{n\geq 1} \chi(n)/n^s$. (This is another
example of a Dirichlet $L$-function.) Suppose that the Generalized
Riemann Hypothesis is true {\it except\/} that $L(s,\chi)$ has a
single zero $\sigma+i\gamma$ with $\sigma>1/2,\ \gamma\geq 0$, and that
$L(s,\bar\chi)$ has a zero at $\sigma-i\gamma$ as well (zeros always
come in conjugate pairs). The formula analogous to (3) is
$$
\frac {  \#\{\text{primes}\ 5n+a \le x\} - \frac14 \pi(x) } {
x^\sigma/\ln x } \approx -2 \text{Re} \left(  \chi(a) \frac
{\cos(\gamma \log x) + i \sin(\gamma \log x)}{\sigma+i\gamma} \right),
$$
for $1\leq a\leq 4$. If we work this out explicitly, we see that

\smallskip\centerline{$\displaystyle \tfrac12(\sigma^2+\gamma^2) \frac
{  \#\{\text{primes}\ 5n+1 \le x\} - \frac14 \pi(x) } { x^\sigma/\ln x
} \approx -\sigma\cos(\gamma\ln x) - \gamma\sin(\gamma\ln x)$}
\smallskip\centerline{$\displaystyle \tfrac12(\sigma^2+\gamma^2) \frac
{  \#\{\text{primes}\ 5n+2 \le x\} - \frac14 \pi(x) } { x^\sigma/\ln x
} \approx \sigma\sin(\gamma\ln x) - \gamma\cos(\gamma\ln x)$}
\smallskip\centerline{$\displaystyle \tfrac12(\sigma^2+\gamma^2) \frac
{  \#\{\text{primes}\ 5n+3 \le x\} - \frac14 \pi(x) } { x^\sigma/\ln x
} \approx -\sigma\sin(\gamma\ln x) + \gamma\cos(\gamma\ln x)$}
\smallskip\centerline{$\displaystyle \tfrac12(\sigma^2+\gamma^2) \frac
{  \#\{\text{primes}\ 5n+4 \le x\} - \frac14 \pi(x) } { x^\sigma/\ln x
} \approx \sigma\cos(\gamma\ln x) + \gamma\sin(\gamma\ln x)$}
\smallskip

Amazingly, this implies that the configuration
$$
\#\{\text{primes}\ 5n+3 \le x\} < \#\{\text{primes}\ 5n+2 \le x\}
< \#\{\text{primes}\ 5n+4 \le x\} < \#\{\text{primes}\ 5n+1 \le x\}
$$
cannot occur when $x$ is sufficiently large! Why? The first inequality
implies that

\medskip
\noindent $-\sigma\sin(\gamma\ln x) + \gamma\cos(\gamma\ln x)
\lessapprox \sigma\sin(\gamma\ln x) - \gamma\cos(\gamma\ln x)$
\smallskip
\hfill ${}\iff 0 \lessapprox \sigma\sin(\gamma\ln x) -
\gamma\cos(\gamma\ln x),$
\medskip

\noindent where the $\lessapprox$ symbol means that the inequality
holds up to an error that goes to zero as $x$ tends to infinity.
Similarly, the third inequality implies that $\sigma\cos(\gamma\ln
x) + \gamma\sin(\gamma\ln x) \lessapprox 0$. The three-inequality
configuration is therefore equivalent to
$$
0 \lessapprox \sigma\sin(\gamma\ln x) - \gamma\cos(\gamma\ln x)
\lessapprox \sigma\cos(\gamma\ln x) + \gamma\sin(\gamma\ln x)
\lessapprox 0,
$$
which implies both $\sigma\sin(\gamma\ln x) - \gamma\cos(\gamma\ln x)
\approx 0$ and $\sigma\cos(\gamma\ln x) + \gamma\sin(\gamma\ln x)
\approx 0$. However, multiplying by $\sin(\gamma\ln x)$ and
$\cos(\gamma\ln x)$ respectively yields $\sigma\sin^2(\gamma\ln x) -
\gamma\cos(\gamma\ln x)\sin(\gamma\ln x) \approx 0$ and
$\sigma\cos^2(\gamma\ln x) + \gamma\cos(\gamma\ln x)\sin(\gamma\ln x)
\approx 0$, and adding these together gives
$$
\sigma = \sigma\sin^2(\gamma\ln x) + \sigma\cos^2(\gamma\ln x) \approx
0.
$$
But this is ridiculous, as it says that the constant $\sigma>\frac12$
tends to 0 as $x$ goes to infinity.

Ford and Konyagin proved more. For
any three arithmetic progressions modulo $q$, there is a feasible
(but equally unlikely) way to prescribe
Generalized-Riemann-Hypothesis-violating zeros of the appropriate
Dirichlet $L$-functions that would cause one of the six orderings
of the three arithmetic progressions not to occur beyond a certain
$x$-value. Moreover, the zeros in these configurations can be
placed arbitrarily far from the real axis, which implies that
there is no way for us to rule out such a possibility by a finite
computation.

It seems that their method can be extended to show that certain
races between just two arithmetic progressions keep the same
leader from some point onwards provided there is another
technically feasible, but rather unlikely, contradiction to the
Generalized Riemann Hypothesis.

    Other similar constructions appear in Ford and Konyagin's work.
Moving up to $r$-way races with $r\ge4$, they have constructed
feasible configurations of ``violations'' which permit no more
than $r^2$ of the $r!$ possible orderings of $r$ arithmetic
progressions to occur infinitely often.

It thus seems that it will  be very difficult to obtain results
independent of any unproved hypothesis.

\bigskip

\remark{Acknowledgements} The section ``Doing the wave'' was
inspired by Enrico Bombieri's delicious article [Bo].
Thanks to Carter Bays, Kevin Ford, Richard Hudson and Nathan Ng
for making available parts of their works that have not appeared in
print. Thanks
also to Guiliana Davidoff and Michael Guy for preparing the appendices.
\endremark

\head Appendix One. A Mount Holyoke College REU.
\endhead

\centerline{ by {\sl Giuliana Davidoff}}

\medskip
\noindent The research group consisted of students   {\sl Caroline
Osowski, Jennifer vanden Eynden, Yi Wang, and Nancy Wrinkle}.

In my senior class in analytic number theory,  I had been asking
students to collect data on various pertinent questions so as to
conjecture a forthcoming theorem, and then to report their results
to the whole class. Having seen many examples where the theorems
matched the initial experimental evidence, I asked them to collect
data on the Mod 4 Race, and then surprised them with Littlewood's
Theorem. Next we  decided to investigate the Mod 3, 5, 7, and 11
Races, in particular whether Team N, the set of primes that are
not congruent to a square mod $q$, consistently leads over Team S,
the set of primes that are  congruent to a square mod $q$. The
results found were intriguingly different from one modulus to
another, so we decided to study this during the upcoming summer
REU program.

The students began by tracking down anything in the existing
literature on sign changes and posting an appeal on a number
theory web site; within days we had received a generous response
from Andrew Odlyzko who pointed us to relevant data of Robert Rumely,
allowing us to begin our own work in earnest.

We were most fascinated by a 1971 paper of Harold Stark [27], in which
he suggested a method to study prime races between primes from any
two given arithmetic progressions. Back then, there was no general
procedure known to prove that each team took the lead infinitely
often in the race between the primes of the form $qn+a$ and the
primes of the form $qn+b$, where $a$ is a square mod $q$ and $b$
is not. Stark's results thus pertain to the first case not covered
by the prior work of Littlewood [15] and of Knapowski and Tur\'an [25,
26].
As he himself pointed out, it seemed particularly difficult to
show $qn+a$ leads infinitely often, even in the case of primes of
the form $5n+4$ racing against primes of the form $5n+2$.

Stark rephrases the problem in terms of two auxiliary functions
(which are complicated expressions, analogous to (3), given in
terms of zeros of various Dirichlet $L$-functions) and, in doing
so, creates a beautiful setting in which he obtains a theorem from
a numerical calculation. In this manner he was able to show that
$\# \{ \text{primes}\ 5n+4\le x\}
> \# \{ \text{primes}\ 5n+2\le x\}$ for arbitrarily large $x$,
assuming the Generalized Riemann Hypothesis.
In fact he proved this with a weaker assumption than
the Generalized Riemann Hypothesis: One needs only that the
Dirichlet $L$-functions associated to this race have no real zeros
in the interval $(1/2,1)$.

The REU group began by making a detailed numerical study of
Stark's formula (for his Mod 5 Race), in order to appreciate how
his complicated expressions vary with the actual numbers of
primes: we found that there are very good correlations, though a
little less so when we use fewer zeros to approximate the formula
analogousto the right-hand side of~(3).

    Using Stark's results and our own numerical calculations, we were 
able
to prove the by-now expected result in the first cases not treated
by him:

\proclaim {Theorem [23]}  Assume that Dirichlet $L$-functions have no 
real zeros
in the interval $(1/2,1)$.
For $a=1$, $2$, or $4$ (the squares mod $q$) and $b=3$, $5$, or $6$ (the
non-squares mod $q$), there are arbitrarily large values of $x$
for which $\# \{ \text{primes}\ 7n+a\le x\}
> \# \{ \text{primes}\ 7n+b\le x\}$. In fact there exists a
constant $c>0$ such that $\# \{ \text{primes}\ 7n+a\le x\}
> \# \{ \text{primes}\ 7n+b\le x\} > c \sqrt{x}/\log x$.
\endproclaim

The proof of the theorem is based on calculations of Stark's
formula using available tables of zeros of Dirichlet
$L$-functions. We could have perhaps gone on to settle this
question about races, one prime modulus at a time.

However, the question that initially interested us was whether
Team S takes the lead over Team N for arbitrarily large $x$. To
study this we had to appropriately modify Stark's formulas: this
time we found that the formulas only involve the zeros of the
Dirichlet $L$-function
$$ \sum_{n\geq 1} \frac {(n/q)}{n^s} $$
where $(n/q)=1$ if $n$ is a square mod $q$,  $(n/q)=0$ if $q$
divides $n$, and $(n/q)=-1$ if $n$ is not a square mod $q$.

We found the same remarkable correlations between the actual
count of primes and our approximation; had the summer not
ended, we would surely have proved that the lead changes hands
infinitely often in this Mod 7 race.

I was later able to go on and prove [23] that the lead changes hands
infinitely often in the Team S vs.~Team N race, for any modulus
$q$, assuming a weak version of the Generalized Riemann
Hypothesis:  One only needs to assume that any real zero of the
above Dirichlet $L$-function lies  to the left of the least upper
bound of the real parts of the complex zeros (which holds, in
particular,
when the Dirichlet $L$-function has no real zeros).

\head Appendix Two. A University of Georgia VIGRE research group
\endhead

\centerline{by {\sl Michael Guy}.}
\medskip
\noindent The research group consisted of students   {\sl Michael
Beck, Zubeyir Cinkir, Michael Guy, Brian Lawler, Eric Pine, Paul
Pollack, and Charles Pooh}, with postdoctoral mentor {\sl Jim
Solazzo}.

The twin prime problem, that there are infinitely many prime pairs
$p, p+2$, is one of the most famous unsolved problems of classical
number theory.  There are many generalizations, but here we shall
focus on the question of whether there are infinitely many prime
pairs $p, p+2k$ for any positive even integer $2k$. The following
conjecture as to how many such pairs there are up to $x$ is due,
essentially, to Hardy and Littlewood [10].

\proclaim{The Hardy-Littlewood conjecture} Let $k$ be a positive
integer, and let   $\pi_{2k}(x)$ be the number of prime pairs
$p,p+2k$ with $p\leq x$. Then
$$
\pi_{2k}(x)\sim  2C_2   \prod\Sb p|k,\ p>2\endSb \frac{p-1}{p-2}
\cdot \Li_{2}(x) ,
$$
where
$$
C_2=\prod\limits_{p>2}\left(1-\frac{1}{(p-1)^2}\right) \text{\rm
and } \Li_2(x)=\int_{2}^{x}\frac{1}{(\log{x})^2}dx .
$$
\endproclaim

We decided to investigate this problem by collecting and analyzing
data. Our program found pairs of primes using a modified
Eratosthenes' sieve and then counted them. Some of the
initial data is collected in the following table:
$$
\matrix
    \underline{X}  &  \underline{\pi_2(X)}  &   \underline{\pi_4(X)}   &
\underline{\pi_6(X)}   &
\underline{\pi_8(X)} &   \underline{\pi_{10}(X)} \\
10^{3} & 35 & 41 & 74 & 38 & 51 \\
10^{4} &    205 & 203& 411 & 208 & 270 \\
10^{5} & 1,224  & 1,216 & 2,447  &   1,260   &  1,624 \\
10^{6} & 8,169   &  8,144  & 16,386  & 8,242  &   10,934 \\
10^{7}  & 58,980  & 58,622  &  117,207 &  58,595  & 78,211 \\
10^{8}  &  440,312  & 440,258  & 879,908  & 439,908 &  586,811 \\
10^{9} & 3,424,506  &   3,424,680  &
6,849,047  & 3,426,124   &  4,567,691\\
10^{10}   & 27,412,679  &
27,409,999  &  54,818,296  & 27,411,508  &  36,548,839 \\
10^{11}   &   224,376,048  & 224,373,161 &
448,725,003  & 224,365,334  & 299,140,330 \\
10^{12}    & 1,870,585,220 & 1,870,585,459  & 3,741,217,498 &
1,870,580,394   &  2,494,056,601\\
\endmatrix
$$
\botcaption{Table 9} $\pi_{2k}(x)$ is the number of prime pairs
$p, p+2k$ with $p\leq x$.
\endcaption
\medskip

Notice that the counts for $2k=2, 4$ and $8$ are very close.
We graphed the data for each $2k\leq 30$, noticing that
$\pi_{2k}(x)$ and $\pi_{2\ell}(x)$ are close if the prime factors
of $2k$ and of $2\ell$ are the same. Indeed,
Hardy and Littlewood's conjecture predicts that $\pi_{2k}(x)\sim
\pi_{2\ell}(x)$ in these circumstances, and so we felt it would be
interesting to study
this ``twin-prime race''.

If we wish to compare $\pi_{2}(x)$ and $\pi_{6}(x)$ then it makes
sense to ``renormalize'' so that the predictions of Hardy and
Littlewood are the same for both. In other words, the conjecture
predicts that $\frac 12 \pi_6(x)$ should be approximately
$\pi_2(x)$, so it makes sense to compare these two quantities.
More generally, we define
$$
\pi_{2k}'(X)=\pi_{2k}(X)\cdot \prod\Sb p|k,\
p>2\endSb\frac{p-2}{p-1} \text{  for each  } k\geq 1,  \text{  and
} \pi_{HL}(X)=2C_2\cdot \Li_{2}(X).
$$
Hardy and Littlewood's  conjecture predicts that $\pi_{2k}'(X)\sim
\pi_{HL}(X)$ for all $k$, so below we have tabulated their
difference.
$$
\matrix
    \underline{X}  & \underline{\pi_{HL}(X)}   &
\underline{\pi_2'-\pi_{HL}}
    &  \underline{\pi_4'-\pi_{HL}} &  \underline{\pi_6'-\pi_{HL}}
    & \underline{\pi_8'-\pi_{HL}}  &\underline{\pi_{10}'-\pi_{HL}} \\
10^{3}  &   45  &  -10  & -4  &  -8 &   -7   & -7\\
10^{4}  & 214  & -9  &  -11 &  -9 &   -6  &  -12\\
10^{5}  & 1,248  &   -24 &  -32 &  -25 &  12  &  -30\\
10^{6}  &   8,248  &   -79 &  -104    &  -55  & -6  &  -48\\
10^{7}  &  58,753 &   227  & -131  &    -150 &     -158   &   -95\\
10^{8}  & 440,367 &  -55  & -109  &    -413 &    -459    &  -259\\
10^{9}  &   3,425,308  &   -802   &   -628  &    -785   &   816  & 460\\
10^{10}  &  27,411,416  &  1,263  &   -1,417 &   -2,268 &   92  &  213\\
10^{11}  & 224,368,866 &  7,182  &   4,295  &   -6,365 & -3,532 &
-13,619\\
10^{12} &    1,870,559,881  &   25,339  &  25,578  &  48,868  &
20,513 & -17,430\\
\endmatrix
$$
\botcaption{Table 10} The renormalized twin prime race.
\endcaption
\medskip

    What a remarkable fit!  It looks like $|\pi_{2k}'(x)-\pi_{HL}(x)|$ is
usually well
less than $\sqrt{x}$. In fact we collected data for
$k=1,2,\ldots,50$ in this range of $x$, and our prediction seemed
to always be very good.

However, this is a paper about prime races after all, and we
wanted to further investigate whether there are any particular
winners, or losers, to this race.  We instructed our program to
take note of when there was a change in first or last position and
to give us this information as well.

\par At the beginning of our investigation, we seemed to have spotted
what we thought were
winners and losers in this prime race.  It appeared that the pairs
$p,p+60$ and $p,p+80$ were ahead of the other pairs with the same
prime divisors, based on data up to $5\cdot 10^9$.  However, after
counting up to $10^{12}$, they were no longer regularly the
winners.  Similarly several prime pairs were consistently the
losers at the start of the prime twin race, but after a while
there were no consistent losers.  As an interesting tidbit, our
program noted that there were 14,455 changes in first place
between $10^3$ and $10^{9}$!

\par So in the end, this appears to be a
race in which there are no particular winners or losers, and still
lots of unanswered questions.

\end